\documentclass[11pt,reqno]{amsart}

\usepackage[
left=0.5in,
right=0.5in,
top=0.7in,
bottom=0.7in
]{geometry}

\usepackage{amsmath,amsthm,amssymb,amscd,amsfonts,mathrsfs,color}
\usepackage[mathscr]{eucal}
\usepackage{mathtools}
\usepackage{pgfplots}
\usepackage{placeins}
\mathtoolsset{showonlyrefs}
\usepgfplotslibrary{fillbetween}
\pgfplotsset{compat=1.18}
\usepackage[bookmarksnumbered,colorlinks,plainpages]{hyperref}
\newtheorem{theorem}{Theorem}[section]
\newtheorem{proposition}[theorem]{Proposition}
\newtheorem{lemma}[theorem]{Lemma}
\newtheorem{corollary}[theorem]{Corollary}
\newtheorem{remark}[theorem]{Remark}
\newtheorem{example}[theorem]{Example}
\theoremstyle{definition}
\newtheorem{definition}[theorem]{Definition}
\numberwithin{equation}{section}

\newcommand{\BH}{\mathcal B(\mathcal H)}
\newcommand{\Real}{\operatorname{Re}}
\newcommand{\Imag}{\operatorname{Im}}
\newcommand{\norm}[1]{\left\lVert #1\right\rVert}
\newcommand{\abs}[1]{\left\lvert #1\right\rvert}
\newcommand{\inner}[2]{\left\langle #1,#2\right\rangle}
\newcommand{\C}{\mathbb C}
\newcommand{\R}{\mathbb R}
\newcommand{\D}{\mathbb D}
\newcommand{\T}{\mathbb T}
\newcommand{\Z}{\mathbb Z}
\newcommand{\conv}{\operatorname{conv}}
\newcommand{\cconv}{\overline{\operatorname{conv}}}

\newcommand{\DW}{\operatorname{DW}}
\newcommand{\dist}{\operatorname{dist}}
\newcommand{\Hol}{\operatorname{Hol}}
\newcommand{\Fab}{\mathcal F_\delta}
\newcommand{\Sch}{\mathcal S_s}
\newcommand{\Tear}{\mathcal T_{\delta,s}}
\newcommand{\Mfr}{\mathfrak M_{\delta,s}}

\begin{document}

\title[Sharp Faber functional calculus on ellipses]
{Sharp Faber functional calculus on ellipses}
\author[M. R. Lahwate]{Mohamed Rassan Lahwate}

\address{Sciences and Technologies Team (ESTE),
			Higher School of Education and Training of El Jadida,
			Chouaib Doukkali University,
			El Jadida, Morocco.}
		\email{mohamedghassane94@gmail.com}


\begin{abstract}
Let $T$ be a Hilbert space operator whose numerical range lies in a closed
elliptical disk $K$. For $g$ in the closed unit ball of the disk algebra, let
$f_g$ be its Faber transform on $K$. We determine the exact supremum of
$\norm{f_g(T)+tI}$ for each prescribed $g(0)$ and each complex shift $t$.
The supremum is attained either at a scalar on the boundary of $K$ or at
a $2\times2$ matrix whose eigenvalues are the foci. Both cases occur.
In particular, the Beckermann--Reichel bound is sharp for every $g(0)$.
The same bound holds for Schur symbols when the spectrum of $T$ lies in
the open ellipse. We obtain the exact union of the numerical ranges and
show that every point of the upper Davis--Wielandt envelope is attained
in dimension three. The shifted norm also gives the sharp
$\rho$-operator-radius bound for $1\leq\rho\leq2$ and a criterion for equality.
\end{abstract}

\maketitle

\section{Introduction}\label{sec:intro}

Let $\mathcal H$ be a complex Hilbert space, and let $\BH$ be the algebra
of bounded linear operators on $\mathcal H$. For $T\in\BH$, the numerical
range and the numerical radius are
\[
 W(T)=\{\inner{Tx}{x}:\norm{x}=1\},
 \qquad
 w(T)=\sup\{\abs{z}:z\in W(T)\};
\]
see~\cite{GustafsonRao}. By von Neumann's inequality~\cite{vonNeumann},
every contraction $T$ satisfies $\norm{g(T)}\leq1$ for every $g$ in the
closed unit ball of the disk algebra $A(\overline\D)$. When the condition
$\norm{T}\leq1$ is replaced by the numerical-range condition
$W(T)\subseteq\overline\D$, the sharp bounds depend on $g(0)$. If
$g(0)=0$, then $w(g(T))\leq1$ by the mapping theorem of Berger and
Stampfli~\cite{BergerStampfli}; see also~\cite{Kato}. For general $g$,
Drury~\cite{Drury} proved the sharp inclusion
\begin{equation}\label{eq:drury}
 W(g(T))
 \subseteq
 \conv\bigl(\overline\D\cup\overline D(g(0),1-\abs{g(0)}^2)\bigr),
\end{equation}
where $\overline D(a,\varrho)$ is the closed disk with center $a$ and
radius $\varrho$, together with the corresponding bound for
$\norm{g(T)}$. The sharpness of this norm bound also follows from the
case $\rho=2$ of the theorem of Schwenninger and de
Vries~\cite{SchwenningerDeVries} recalled below.

Klaja et al.~\cite{KlajaMashreghiRansford} gave a short proof of the
inclusion~\eqref{eq:drury} based on the local estimate
\[
 \norm{Tx}^2\leq2+2\sqrt{1-\abs{\inner{Tx}{x}}^2},
 \qquad
 \norm{x}=1,
\]
that is, from the geometry of the Davis--Wielandt shell
\[
 \DW(T)=\{(\inner{Tx}{x},\norm{Tx}^2):\norm{x}=1\}\subset\C\times\R;
\]
see~\cite{Wielandt,Davis}. The same two moments control the operator
radii $w_\rho$ of Sz.-Nagy and Foias~\cite{SzNagyFoias}, which
interpolate between the norm $(\rho=1)$ and the numerical radius
$(\rho=2)$; see Subsection~\ref{ssec:radii}. Through the deformed
numerical ranges of Pagacz et al.~\cite{PagaczPietrzyckiWojtylak}, the
shell determines $w_\rho(T)$ as an outer radius. On the disk, the sharp
$g(0)$-dependent norm bound for $\rho$-contractions is due to
Schwenninger and de Vries~\cite{SchwenningerDeVries}, with a
matrix-valued extension by Hartz and de Vries~\cite{HartzDeVries}.

We determine the corresponding sharp calculus when the disk is replaced
by an ellipse. For $0\leq\delta<1$, let $K_\delta$ be the closed
elliptical disk with semi-axes $1\pm\delta$ and foci $\pm2\sqrt\delta$,
and let $G_\delta$ be its interior. Thus $K_0=\overline\D$, and every
nondegenerate closed elliptical disk is the image of some $K_\delta$
under an affine change of variable. The functional calculus of the class
$W(T)\subseteq K_\delta$ is measured by the B.~and F.~Delyon
norm~\cite{DelyonDelyon}
\[
 \norm{f}_{\mathrm{bfd}}
 =
 \sup\bigl\{\norm{f(T)}:\overline{W(T)}\subset G_\delta\bigr\},
\]
which is at most $(1+\sqrt2)\sup_{G_\delta}\abs{f}$ by the theorem of
Crouzeix and Palencia~\cite{CrouzeixPalencia}; the sharp constant $2$
conjectured by Crouzeix~\cite{Crouzeix} has recently been
announced~\cite{Jin,LoristSchwenninger}. Agler
et al.~\cite{AglerLykovaYoung,AglerLykovaYoungBfd} related this norm to
a functional calculus on an annulus and on the symmetrized bidisk. Their
key step is the annular representation (Theorem~\ref{thm:aly}): if
$0<\delta<1$ and $T$ is a matrix with $W(T)\subseteq K_\delta$, then
\begin{equation}\label{eq:aly-intro}
 T=(X+\delta X^{-1})\big|_{\mathcal H},
 \qquad
 \norm{X}\leq1,
 \qquad
 \norm{\delta X^{-1}}\leq1,
\end{equation}
for an invertible operator $X$ on $\mathcal H\oplus\mathcal H$ such that
$X+\delta X^{-1}$ leaves $\mathcal H\cong\mathcal H\oplus\{0\}$
invariant. Thus $T$ is a part of $\Psi_\delta(X)$, where
$\Psi_\delta(w)=w+\delta/w$ is the Joukowski map.

This representation explains why Faber polynomials appear. Since
$\Psi_\delta(\delta/w)=\Psi_\delta(w)$, the pullback of a function $f$
holomorphic on $G_\delta$ to the annulus $\{\delta<\abs{w}<1\}$ is
invariant under $w\mapsto\delta/w$. Its Laurent expansion therefore
groups into the pairs $w^n+(\delta/w)^n$, which are the pullbacks of the
Faber polynomials $F_{n,\delta}$ of $K_\delta$~\cite{Faber1903,Faber1920}.
Hence $f=f_g$, where
\begin{equation}\label{eq:faber-intro}
 f_g\Bigl(w+\frac{\delta}{w}\Bigr)
 =
 g(w)+g\Bigl(\frac{\delta}{w}\Bigr)-g(0),
 \qquad
 \delta<\abs{w}<1,
\end{equation}
for a unique function $g$ holomorphic on $\D$
(Remark~\ref{rem:laurent}). We call $f_g$ the Faber transform of $g$;
for $\delta=0$, one has $f_g=g$. For $\delta>0$, the value $g(0)$ is the
constant Faber coefficient of $f_g$, not $f_g(0)$; for instance, the
symbols $u$ and $u^2$ have transforms $z$ and $z^2-2\delta$.

If $g$ is holomorphic on a neighborhood of $\overline\D$,
then~\eqref{eq:aly-intro} and~\eqref{eq:faber-intro} give
\[
 f_g(T)=\bigl(g(X)+g(\delta X^{-1})-g(0)I\bigr)\big|_{\mathcal H}.
\]
Thus the elliptic problem becomes a problem for the two commuting
contractions $X$ and $\delta X^{-1}$, coupled through the single value
$g(0)$. This leads to the symbol class
\begin{equation}\label{eq:schur-class}
 \Sch=\{g\in H^\infty(\D):\norm{g}_\infty\leq1,\ g(0)=s\},
 \qquad
 0\leq s\leq1.
\end{equation}

Uniform estimates for the Faber calculus are known for every convex set
$E$. Beckermann~\cite{Beckermann} proved that $\norm{F_n(T)}\leq2$ for
the Faber polynomials $F_n$ of $E$ whenever $W(T)\subseteq E$, extending
a scalar bound of K\H{o}v\'ari and Pommerenke~\cite{KovariPommerenke}.
Beckermann and Reichel~\cite[Theorem~2.1]{BeckermannReichel} proved, for
disk-algebra symbols,
\begin{equation}\label{eq:br}
 \norm{f_g(T)+g(0)I}\leq2\norm{g}_\infty.
\end{equation}
These estimates normalize the symbol, as we do, but they do not depend
on $g(0)$, and they do not locate $W(f_g(T))$. Putinar and
Tsikalas~\cite{PutinarTsikalas} obtained complete bounds, also on annuli
and through Faber $\rho$-dilations, with ellipses among their examples;
their estimates control the calculus through the supremum of the
function on its domain. In contrast, we prescribe the norm and the
constant coefficient of the disk symbol and compute the extremum for
every complex shift.

One family of norms answers both questions. For $F=f_g(T)$ and a unit
vector $x$,
\[
 \norm{(F+tI)x}^2
 =
 \norm{Fx}^2+2\Real\bigl(\overline t\inner{Fx}{x}\bigr)+\abs{t}^2,
 \qquad
 t\in\C.
\]
As $\abs{t}\to\infty$, the shifted norms recover the support function
of $W(F)$. At finite $t$, they bound the second moment $\norm{Fx}^2$
against the first, that is, they describe the upper boundary of the
Davis--Wielandt shell. At $t=g(0)$, they give the quantity
in~\eqref{eq:br}. After a rotation of $g$ and $t$, we may take
$g(0)=s\geq0$. We therefore study, for every $t\in\C$, the shifted norm
\begin{equation}\label{eq:shifted-intro}
 \Mfr(t)=\sup_{(g,T)}\norm{f_g(T)+tI}
\end{equation}
over the admissible pairs: $g\in\Sch$, $W(T)\subseteq K_\delta$, and
either $g\in A(\overline\D)$ or $\sigma(T)\subset G_\delta$
(Definition~\ref{def:admissible}). For $s=1$ the only symbol is the
constant $1$, so that $\Mfr(t)=\abs{1+t}$. For $s<1$, put $k=s^2$. Two
models suffice: a boundary scalar and the focal matrix
\begin{equation}\label{eq:focal-matrix}
 A_\delta=\begin{pmatrix}0&2\\2\delta&0\end{pmatrix},
\end{equation}
whose numerical range is $K_\delta$ and whose eigenvalues are the foci
$\pm2\sqrt\delta$.

\begin{theorem}\label{thm:main}
Let $0\leq\delta<1$ and $0\leq s<1$. For every $t\in\C$,
\begin{equation}\label{eq:main-intro}
 \Mfr(t)
 =
 \max\Bigl\{R+\abs{C+t},\ r+\sqrt{d^2+\abs{c+t}^2}\Bigr\},
\end{equation}
where
\begin{equation}\label{eq:parameters}
 C=-s\delta b,
 \qquad
 R=1+b,
 \qquad\qquad
 c=s-2s\delta a,
 \qquad
 r=(1+\delta)a,
 \qquad
 d=(1-\delta)a,
\end{equation}
with
\[
 a=\frac{1-k}{1-k\delta},
 \qquad
 b=\frac{(1-k)\delta}{1-k\delta^2}.
\]
The first value is attained by a scalar on $\partial K_\delta$ and the
second by $A_\delta$. In each case the symbol can be chosen from the
automorphisms $g_\eta(u)=(s+\eta u)/(1+s\eta u)$, $\eta\in\T$.
The phase may depend on the model and on $t$.
\end{theorem}

Even at $t=0$, either model can give the larger norm. For
$(\delta,s)=(0.1,0.5)$ the focal value is $1.658$ and the scalar value
is $1.079$, to three decimal places; for $(\delta,s)=(0.9,0.9)$ the
scalar value is $1.900$ and the focal value is $1.578$. At $t=s$, the
focal value equals $2$ for all $\delta$ and $s$. Thus
Theorem~\ref{thm:main} gives $\Mfr(s)=2$: the Beckermann--Reichel
bound~\eqref{eq:br} is sharp for every prescribed $g(0)$, and the
theorem evaluates the supremum at every other shift.

The proof of the upper bound combines~\eqref{eq:aly-intro} with the
quantum annulus dilation of McCullough and Pascoe~\cite{McCulloughPascoe}.
This dilates $f_g(T)$ to $g(U)+g(\delta U^{-1})-sI$, where $\abs{U}$
takes only the values $\delta$ and $1$. A Heinz-type interpolation and a
Schur complement then give a block inequality whose entries are
functions of the single self-adjoint operator
$(\abs{U^*}-\abs{U})/(1+\delta)$. By the spectral theorem, the estimate
reduces to a $2\times2$ scalar matrix inequality on an interval. After a
positive denominator is cleared, the determinant of this matrix is
affine in the spectral variable, so only the two endpoints of the
interval matter. They give exactly the scalar and the focal values. The
same argument gives the sharp bound for $g(X)+g(\delta X^{-1})$ directly
on the quantum annulus (Corollary~\ref{cor:annulus}).

The shifted norms also determine the geometry of the outputs. Letting
$\abs{t}\to\infty$ gives the closed convex hull of the numerical ranges,
and we prove that every point of this hull is attained, in dimension at
most three:
\begin{equation}\label{eq:teardrop-intro}
 \bigcup_{(g,T)}W(f_g(T))
 =
 \Tear:=\conv\bigl(\overline D(C,R)\cup\overline D(c,r)\bigr).
\end{equation}
For $\delta=0$, the two disks have centers $0,s$ and radii $1,1-s^2$,
which gives Drury's teardrop~\eqref{eq:drury}; for $s=0$, the teardrop
is the disk of radius $1+\delta$. Thus, in contrast with the
Berger--Stampfli theorem, $W(f_g(T))$ can leave $K_\delta$ when
$\delta>0$, even if $g(0)=0$ (Remark~\ref{rem:teardrop}). At finite shifts we obtain the upper
envelope of the Davis--Wielandt shells (Theorem~\ref{thm:dw}). Every
point on it is realized in dimension three by one symbol $g_\eta$ and an
operator $A_\delta\oplus\lambda$ with $\lambda\in\partial K_\delta$. On
a common supporting plane, the same phase $\eta$ realizes both contacts.

For $1\leq\rho\leq2$, put $p=\rho-1$ and $q=2-\rho$. With the
parameters of Theorem~\ref{thm:main}, the sharp output $\rho$-radius over
the same admissible pairs is
\[
 \sup_{(g,T)}w_\rho(f_g(T))
 =
 \max\left\{R-C,\
 \frac{r+p\abs{c}+\sqrt{(r+p\abs{c})^2-\rho q(r^2-d^2-c^2)}}{\rho}
 \right\}.
\]
For the single symbol $g_1(u)=(s+u)/(1+su)$, the boundary scalar
$-(1+\delta)$ and the focal matrix $A_\delta$ realize the two values,
respectively (Corollary~\ref{cor:radius}). Here $\rho$ measures only the
output; the input is always the condition $W(T)\subseteq K_\delta$. For
$\delta=0$, the formula gives Drury's bound $1+s-s^2$ for $w(g(T))$ when
$\rho=2$, and the norm bound $1-s^2+\sqrt{(1-s^2)^2+s^2}$ when $\rho=1$;
the latter is the case $T\in\mathcal C_2$ of~\cite{SchwenningerDeVries}.
When the input is instead $w_\rho(T)\leq1$ on the disk, the sharp
$g(0)$-dependent bound for $w_\rho(g(T))$ is obtained
in~\cite{LahwateIghachane}. The two problems overlap only for
$\delta=0$ and $\rho=2$, where both give Drury's bound.

Finally, we characterize equality. For an admissible pair,
$\norm{f_g(T)+tI}=\Mfr(t)$ exactly when the closed Davis--Wielandt shell
of $f_g(T)$ meets an explicit exposed face of the envelope
(Theorem~\ref{thm:equality-contact}). For a fixed linear-fractional
transform $f$ in finite dimension, equality in the $\rho$-radius bound
forces a reducing boundary scalar when the scalar value is strictly
larger, and a reducing copy of $A_\delta$ when the focal value is
strictly larger (Theorem~\ref{thm:lf-equality}). This rigidity is
special to linear-fractional transforms: for the symbol $u^2$ on the
disk, a $3\times3$ weighted shift attains the strict focal value without
a reducing copy of $A_0$ (Remark~\ref{rem:not-lf}).

Section~\ref{sec:pre} collects the preliminaries and the annular
dilation of the Faber calculus. Section~\ref{sec:block} proves
Theorem~\ref{thm:main}: the block inequality gives the upper bound, and
the scalar and focal models attain it. Section~\ref{sec:consequences}
derives the elliptic teardrop, the upper Davis--Wielandt envelope and the
sharp operator radii. Section~\ref{sec:equality} characterizes equality,
first by contact and then for linear-fractional transforms.

\section{Faber calculus on an ellipse and annular dilation}\label{sec:pre}

Throughout the paper, $0\leq\delta<1$, $0\leq s\leq1$ and $k=s^2$.
Whenever operator radii occur, $1\leq\rho\leq2$, and we write
\[
 p=\rho-1,
 \qquad
 q=2-\rho;
 \qquad
 \text{thus}\quad p,q\geq0
 \quad\text{and}\quad
 p^2+\rho q=1.
\]
The inner product of $\mathcal H$ is linear in the first variable. We
write $\D$ for the open unit disk, $\T$ for the unit circle, and
$\Hol(\overline\D)$ for the set of functions holomorphic on some disk
$\{\abs{u}<R_0\}$ with $R_0>1$. For a subset $S$ of a real vector space,
$\conv S$ is its convex hull and $\cconv S$ its closed convex hull.

\subsection{Numerical ranges and operator radii}\label{ssec:radii}

For $T\in\BH$, the numerical-range resolvent estimate
$\norm{(zI-T)^{-1}}\leq1/\dist(z,W(T))$ holds outside
$\overline{W(T)}$, and $\sigma(T)\subseteq\overline{W(T)}$;
see~\cite{GustafsonRao}. Hence a function holomorphic near $K_\delta$
has a holomorphic functional calculus at every $T$ with
$W(T)\subseteq K_\delta$.

We also use that the holomorphic functional calculus is continuous. If
$f_j\to f$ locally uniformly on a neighborhood of $\sigma(S)$, then
$f_j(S)\to f(S)$. Moreover, $S\mapsto f(S)$ is norm continuous on the
set of operators whose spectrum lies in the domain of $f$;
see~\cite[Section~VII.4]{Conway}.

For $\rho>0$, the class $\mathcal C_\rho$ of Sz.-Nagy and
Foias~\cite{SzNagyFoias} consists of operators $T$ admitting a unitary
$\widetilde U$ on a space $\mathcal K\supseteq\mathcal H$ such that
\begin{equation}\label{eq:rho-dilation}
 T^n=\rho P_{\mathcal H}\widetilde U^n\big|_{\mathcal H},
 \qquad
 n\geq1.
\end{equation}
The associated operator radius is
\begin{equation}\label{eq:rho-radius}
 w_\rho(T)=\inf\{r>0:r^{-1}T\in\mathcal C_\rho\}.
\end{equation}
We use the following positivity criterion for $1\leq\rho\leq2$.

\begin{lemma}\label{lem:pencil}
Let $1\leq\rho\leq2$ and $S\in\BH$.

\noindent{\rm (i)} For $h>0$, one has $w_\rho(S)\leq h$ if and only if
\begin{equation}\label{eq:pencil}
 \rho h^2I-2ph\Real(\overline\zeta S)-qS^*S\geq0,
 \qquad
 \zeta\in\T.
\end{equation}

\noindent{\rm (ii)} The function $w_\rho$ is a norm on $\BH$, with
$w_1(S)=\norm{S}$, $w_2(S)=w(S)$ and
$w(S)\leq w_\rho(S)\leq\norm{S}$. Moreover, $w_\rho(\lambda I)=\abs{\lambda}$
and $w_\rho(S_1\oplus S_2)=\max\{w_\rho(S_1),w_\rho(S_2)\}$.
\end{lemma}

\begin{proof}
By~\cite[Section~I.11]{SzNagyFoiasBook}, applied to $h^{-1}S$, one has
$w_\rho(S)\leq h$ if and only if
\[
 \rho h^2I-2ph\Real(\overline zS)-q\abs{z}^2S^*S\geq0,
 \qquad
 z\in\overline\D.
\]
For fixed $\zeta\in\T$ and a vector $x$, the quadratic form of this
operator at $z=\varrho\zeta$ is a concave function of $\varrho\in[0,1]$,
because $q\geq0$. It equals $\rho h^2\norm{x}^2\geq0$ at $\varrho=0$.
Hence it is nonnegative on $[0,1]$ as soon as it is nonnegative at
$\varrho=1$. This proves part~{\rm (i)}. The first assertions of
part~{\rm (ii)} are proved in~\cite{Holbrook,AndoNishio}.

The scalar identity follows from $w\leq w_\rho\leq\norm{\cdot}$.
The pencil in~\eqref{eq:pencil} respects direct sums, which gives
the last assertion.
\end{proof}

Recall that the Davis--Wielandt shell records the two moments that occur
in the quadratic form of~\eqref{eq:pencil}:
\[
 \DW(S)=\{(\inner{Sx}{x},\norm{Sx}^2):\norm{x}=1\}.
\]
Pagacz et al.~\cite{PagaczPietrzyckiWojtylak} use these moments to define a deformed numerical range. Define
$\Phi_\rho:\C\times[0,\infty)\to\C$ by
\begin{equation}\label{eq:phi}
 \Phi_\rho(z,b)
 =
 \begin{cases}
 \displaystyle
 \frac{z}{\rho\abs{z}}
 \Bigl(p\abs{z}+\sqrt{p^2\abs{z}^2+\rho qb}\Bigr),
 & z\neq0,\\[6pt]
 0,
 & z=0,
 \end{cases}
\end{equation}
and put
\[
 W^\rho(S)=\cconv\Phi_\rho(\DW(S)).
\]
For $\rho=2$, the map $\Phi_2(z,b)=z$ forgets the second moment. For
$\rho=1$ and $z\neq0$, the map $\Phi_1(z,b)=\sqrt{b}\,z/\abs{z}$ keeps only the
direction of $z$ and the length $\sqrt{b}=\norm{Sx}$.
By~\cite[Theorem~9 and Corollary~10]{PagaczPietrzyckiWojtylak},
\begin{equation}\label{eq:ppw}
 W^2(S)=\overline{W(S)},
 \qquad
 w_\rho(S)=\max\{\abs{z}:z\in W^\rho(S)\}.
\end{equation}

Two elementary properties of $\Phi_\rho$ will be used. Since $p+1=\rho$
and $p^2+\rho q=1$, one has
\[
 \Phi_\rho(z,\abs{z}^2)=z.
\]
For fixed $z\neq0$, the point $\Phi_\rho(z,b)$ moves outward along the
ray from $0$ through $z$ as $b$ increases.

\subsection{Elliptical domains and the focal matrix}\label{ssec:ellipse}

The closed elliptical disk $K_\delta$ is described by its gauge function
\[
 m_\delta(x+iy)
 =
 \Bigl(\frac{x^2}{(1+\delta)^2}+\frac{y^2}{(1-\delta)^2}\Bigr)^{1/2},
\]
so that
\[
 K_\delta=\{m_\delta\leq1\},
 \qquad
 G_\delta=\{m_\delta<1\}.
\]
The foci are $\pm2\sqrt\delta$ because $(1+\delta)^2-(1-\delta)^2=4\delta$.

The real-linear map
\[
 L_\delta(\chi)=\chi+\delta\overline\chi,
 \qquad
 \chi\in\C,
\]
sends $\varrho e^{i\theta}$ to
\[
 \varrho(1+\delta)\cos\theta+i\varrho(1-\delta)\sin\theta.
\]
Thus $L_\delta$ is a bijection of $\C$, with inverse
\[
 L_\delta^{-1}(z)=\frac{z-\delta\overline z}{1-\delta^2}.
\]
It satisfies $m_\delta(L_\delta(\chi))=\abs{\chi}$, and it maps
$\overline\D$ onto $K_\delta$.

We next consider the focal matrix~\eqref{eq:focal-matrix}. For
$\delta=0$, it is twice the $2\times2$
Jordan block, the classical extremal operator for the condition
$w(T)\leq1$; see~\cite{Ando}. In general, its numerical range and its second moments are
computed through $L_\delta$.

\begin{lemma}\label{lem:focal}
For a unit vector $x=(x_1,x_2)\in\C^2$, put $\chi=2x_2\overline{x_1}$.
Then
\[
 \inner{A_\delta x}{x}=L_\delta(\chi),
 \qquad
 \norm{A_\delta x}^2=4\abs{x_2}^2+4\delta^2\abs{x_1}^2.
\]
Moreover, $W(A_\delta)=K_\delta$ and $A_\delta^2=4\delta I$. For
$\delta>0$, the eigenvalues of $A_\delta$ are the foci $\pm2\sqrt\delta$.
\end{lemma}

\begin{proof}
The two moment identities follow from $A_\delta x=(2x_2,2\delta x_1)$.
For $x=(\cos\alpha,e^{i\theta}\sin\alpha)$ the parameter
$\chi=e^{i\theta}\sin 2\alpha$ ranges over $\overline\D$ as
$0\leq\alpha\leq\pi/4$ and $\theta$ varies. Conversely,
$\abs{\chi}\leq\abs{x_1}^2+\abs{x_2}^2=1$ for every unit vector.
Thus $W(A_\delta)=L_\delta(\overline\D)=K_\delta$.
Matrix multiplication gives $A_\delta^2=4\delta I$.
\end{proof}

For $f$ holomorphic on $G_\delta$, put
\begin{equation}\label{eq:bfd}
 \norm{f}_{\rho,\mathrm{bfd}}
 =
 \sup\{w_\rho(f(T)):\overline{W(T)}\subset G_\delta\}.
\end{equation}
At $\rho=1$ this is the B.~and F.~Delyon norm~\cite{DelyonDelyon},
studied on the ellipse in~\cite{AglerLykovaYoung,AglerLykovaYoungBfd}.
When $f$ extends across the boundary, the supremum also allows the
boundary scalars and the focal matrix.

\begin{lemma}\label{lem:closed-class}
If $f$ is holomorphic on a neighborhood of $K_\delta$, then
\[
 \norm{f}_{\rho,\mathrm{bfd}}
 =
 \sup\bigl\{w_\rho(f(T)):W(T)\subseteq K_\delta\bigr\}.
\]
\end{lemma}

\begin{proof}
The inequality ``$\leq$'' holds because $\overline{W(T)}\subset G_\delta$
implies $W(T)\subseteq K_\delta$.

Conversely, let $W(T)\subseteq K_\delta$. For $0<\theta<1$, one has
\[
 \overline{W(\theta T)}\subseteq\theta K_\delta\subset G_\delta.
\]
Moreover, $f(\theta T)\to f(T)$ in norm as $\theta\uparrow1$, by the continuity of
the functional calculus. Since $w_\rho$ is a norm with
$w_\rho\leq\norm{\cdot}$, it follows that $w_\rho(f(\theta T))\to w_\rho(f(T))$.
\end{proof}

\subsection{Faber polynomials and the Faber transform}\label{ssec:faber}

For $0<\delta<1$, the Joukowski map $\Psi_\delta(w)=w+\delta/w$
sends $\{\abs{w}>1\}$
conformally onto $\C\setminus K_\delta$, and
$\Psi_\delta(w)=w+O(1/w)$ at infinity. For $\varrho\geq\sqrt\delta$, it maps the
circle $\{\abs{w}=\varrho\}$ onto the ellipse with semi-axes $\varrho+\delta/\varrho$ and
$\varrho-\delta/\varrho$ and foci $\pm2\sqrt\delta$; for $\varrho=1$, this ellipse is
$\partial K_\delta$. For $\delta=0$, we put $\Psi_0(w)=w$.

The Faber polynomials of $K_\delta$ are defined by the generating
function~\cite{Faber1903,Curtiss}
\begin{equation}\label{eq:faber-gen}
 \frac{\Psi_\delta'(w)}{\Psi_\delta(w)-z}
 =
 \sum_{n=0}^\infty\frac{F_{n,\delta}(z)}{w^{n+1}},
 \qquad \abs{w}\ \text{large}.
\end{equation}
To compute them, fix $z\in\C$ and let $u_1(z),u_2(z)$ be the roots of
\begin{equation}\label{eq:roots}
 u^2-zu+\delta=0,
 \qquad\text{so that}\quad
 u_1+u_2=z,
 \quad
 u_1u_2=\delta,
\end{equation}
labeled so that $\abs{u_1}\geq\abs{u_2}$. Since
\[
 \Psi_\delta(w)-z=\frac{(w-u_1)(w-u_2)}{w},
 \qquad
 \Psi_\delta'(w)=\frac{w^2-\delta}{w^2},
\]
one has
\begin{align*}
 \frac{\Psi_\delta'(w)}{\Psi_\delta(w)-z}
 &=
 \frac{w^2-\delta}{w(w-u_1)(w-u_2)}\\
 &=
 \frac{1}{w-u_1}+\frac{1}{w-u_2}-\frac{1}{w}.
\end{align*}
Expanding in powers of $1/w$ and comparing with~\eqref{eq:faber-gen}, we
obtain $F_{0,\delta}=1$ and
\begin{equation}\label{eq:faber-poly}
 F_{n,\delta}(z)=u_1(z)^n+u_2(z)^n,
 \qquad
 F_{n,\delta}\Bigl(w+\frac{\delta}{w}\Bigr)=w^n+\Bigl(\frac{\delta}{w}\Bigr)^n,
 \qquad
 n\geq1.
\end{equation}

Thus $F_{n,\delta}$ is the $n$-th power sum of the two roots.
Equivalently,
\[
 F_{1,\delta}(z)=z,
 \qquad
 F_{2,\delta}(z)=z^2-2\delta,
 \qquad
 F_{n+1,\delta}=zF_{n,\delta}-\delta F_{n-1,\delta},
 \quad n\geq2.
\]
For $\delta>0$, one has $F_{n,\delta}(z)=2\delta^{n/2}T_n\bigl(z/(2\sqrt\delta)\bigr)$,
where $T_n$ is the Chebyshev polynomial of the first
kind~\cite{Faber1920,Curtiss,Suetin}. For $\delta=0$, the roots are $z$
and $0$, and~\eqref{eq:faber-poly} gives $F_{n,0}(z)=z^n$.

The position of the roots is measured by the function
\[
 \varrho_\delta(z)=\abs{u_1(z)}.
\]

\begin{lemma}\label{lem:levels}
The function $\varrho_\delta$ is continuous, $\abs{u_2}\leq\varrho_\delta$,
and
\[
 K_\delta=\{\varrho_\delta\leq1\},
 \qquad
 G_\delta=\{\varrho_\delta<1\}.
\]
In particular, $\varrho_\delta$ is bounded by a constant smaller than one
on each compact subset of $G_\delta$.
\end{lemma}

\begin{proof}
For $\delta=0$, one has $\varrho_0(z)=\abs{z}$.

Let $\delta>0$. The unordered pair of roots depends continuously on $z$,
so their larger modulus $\varrho_\delta(z)$ is continuous. Moreover,
\[
 \abs{u_1}^2\geq\abs{u_1u_2}=\delta.
\]
Write $u_1=\varrho e^{i\theta}$ with $\varrho=\varrho_\delta(z)\geq\sqrt\delta$. Then
\begin{align*}
 z
 &=
 u_1+\frac{\delta}{u_1}\\
 &=
 \Bigl(\varrho+\frac{\delta}{\varrho}\Bigr)\cos\theta
 +i\Bigl(\varrho-\frac{\delta}{\varrho}\Bigr)\sin\theta,
\end{align*}
so $z$ lies on the ellipse $\Psi_\delta(\varrho\T)$. For $\varrho=\sqrt\delta$, this
ellipse degenerates to the focal segment $[-2\sqrt\delta,2\sqrt\delta]$.

Both semi-axes $\varrho\pm\delta/\varrho$ increase strictly with $\varrho\geq\sqrt\delta$.
These confocal ellipses are therefore nested and pairwise disjoint. The closed
region $K_\delta$ is the union of those with $\varrho\leq1$, and $G_\delta$ is the
union of those with $\varrho<1$. Consequently, $z\in K_\delta$ if and only if
$\varrho_\delta(z)\leq1$, and $z\in G_\delta$ if and only if
$\varrho_\delta(z)<1$. The last assertion follows from the continuity of
$\varrho_\delta$.
\end{proof}

By Lemma~\ref{lem:levels}, both roots lie in $\D$ when $z\in G_\delta$.
This makes the following definition possible.

\begin{definition}\label{def:faber}
Let $g(u)=\sum_{n\geq0}a_nu^n$ be holomorphic on $\D$. Its \emph{Faber
transform} with respect to $K_\delta$ is the function
\[
 f_g(z)
 =
 (\Fab g)(z)
 =
 g(u_1(z))+g(u_2(z))-g(0),
 \qquad
 z\in G_\delta.
\]
\end{definition}

The right-hand side is symmetric in $u_1$ and $u_2$, so it does not
depend on the labeling of the roots. The next lemma shows that $f_g$ is
the function of~\eqref{eq:faber-intro}.

\begin{lemma}\label{lem:faber-transform}
Let $g(u)=\sum_na_nu^n$ be holomorphic on $\D$.

\noindent{\rm (i)} The function $f_g$ is holomorphic on $G_\delta$, the
series $a_0+\sum_{n\geq1}a_nF_{n,\delta}$ converges to $f_g$ locally
uniformly on $G_\delta$, and
\[
 f_g\Bigl(w+\frac{\delta}{w}\Bigr)
 =
 g(w)+g\Bigl(\frac{\delta}{w}\Bigr)-g(0),
 \qquad
 \delta<\abs{w}<1.
\]
For $\delta=0$, one has $f_g=g$. If $g$ is a polynomial of degree at most
$N$, then
\[
 f_g=a_0+\sum_{n=1}^Na_nF_{n,\delta}
\]
is a polynomial of degree at most $N$.

\noindent{\rm (ii)} The map $\Fab$ is linear, $\Fab1=1$, and it
satisfies the a priori bound
\[
 \sup_{G_\delta}\abs{f_g}\leq2\norm{g}_\infty+\abs{g(0)}.
\]
\end{lemma}

\begin{proof}
For part~{\rm (i)}, fix $0<\varrho_0<1$, and choose $\varrho_0<\varrho_1<1$. The Cauchy
estimates give
\[
 \abs{a_n}\leq \varrho_1^{-n}\max_{\abs{u}=\varrho_1}\abs{g},
\]
and hence $\sum_n\abs{a_n}\varrho_0^n<\infty$.

On the set $\{\varrho_\delta\leq\varrho_0\}$, both roots have modulus at most $\varrho_0$.
Therefore~\eqref{eq:faber-poly} gives
\[
 \abs{a_nF_{n,\delta}}\leq2\abs{a_n}\varrho_0^n
\]
there. Consequently, the series $a_0+\sum_{n\geq1}a_nF_{n,\delta}$
converges uniformly on $\{\varrho_\delta\leq\varrho_0\}$, and its sum is
\[
 g(u_1)+g(u_2)-a_0=f_g.
\]
Since the partial sums are polynomials, $f_g$ is holomorphic on the open
set $\{\varrho_\delta<\varrho_0\}$. By Lemma~\ref{lem:levels}, every compact subset
of $G_\delta$ lies in $\{\varrho_\delta\leq\varrho_0\}$ for some $\varrho_0<1$, and the
union of the sets $\{\varrho_\delta<\varrho_0\}$ is $G_\delta$. This proves the
first assertion.

The formula on the annulus holds because the roots of
$u^2-(w+\delta/w)u+\delta$ are $w$ and $\delta/w$. For $\delta=0$, the
roots are $z$ and $0$, and therefore $f_g=g$. If $g$ is a polynomial, the
series is a finite sum, and $F_{n,\delta}$ has degree $n$ by the
recurrence above.

Part~{\rm (ii)} follows from the definition, since
$\abs{g(u_i)}\leq\norm{g}_\infty$.
\end{proof}

\begin{remark}\label{rem:laurent}
Let $0<\delta<1$. Every function $f$ holomorphic on $G_\delta$ is the
Faber transform of a unique $g$ holomorphic on $\D$. Indeed,
$\Psi_\delta$ maps the annulus
\[
 \mathbb A_\delta=\{\delta<\abs{w}<1\}
\]
onto $G_\delta$, and $\Psi_\delta(\delta/w)=\Psi_\delta(w)$. Hence
$h=f\circ\Psi_\delta$ is holomorphic on $\mathbb A_\delta$ and
$h(\delta/w)=h(w)$. If $h(w)=\sum_{n\in\Z}a_nw^n$, this symmetry means
$a_{-n}=\delta^na_n$, so that
\[
 h(w)=a_0+\sum_{n\geq1}a_n\Bigl(w^n+\frac{\delta^n}{w^n}\Bigr).
\]

The series $g(u)=\sum_{n\geq0}a_nu^n$ converges on $\D$, and
Lemma~\ref{lem:faber-transform}{\rm(i)} gives $f=f_g$. Uniqueness holds
because $f$ determines the coefficients $a_n$. In particular,
\[
 g(0)
 =
 a_0
 =
 \frac{1}{2\pi}\int_0^{2\pi}
 f\Bigl(\varrho e^{i\theta}+\frac{\delta}{\varrho}e^{-i\theta}\Bigr)\,d\theta,
 \qquad
 \delta<\varrho<1,
\]
is the mean value of $f$ over each confocal ellipse $\Psi_\delta(\varrho\T)$ in
its Joukowski parametrization.

The normalization $\norm{g}_\infty\leq1$ used below is a condition on the
symbol, not on $f$. For example, $g(u)=u$ gives $f_g(z)=z$, whose
supremum on $G_\delta$ is $1+\delta$.
\end{remark}

\subsection{Annular representation and dilation}\label{sec:dilation}

The annular representation of Agler et al.~\cite{AglerLykovaYoung} uses
an invertible contraction $X$ for which $\delta X^{-1}$ is also a
contraction. These are the Douglas--Paulsen
operators~\cite{DouglasPaulsen}; see~\cite{AglerLykovaYoungDp} for the
associated function theory on the annulus.

\begin{theorem}[{\cite[Proposition~10.4]{AglerLykovaYoung}}]\label{thm:aly}
Let $0<\delta<1$, $\dim\mathcal H<\infty$ and $W(T)\subseteq K_\delta$.
Then there is an invertible operator
$X\in\mathcal B(\mathcal H\oplus\mathcal H)$ with $\norm{X}\leq1$ and
$\norm{X^{-1}}\leq1/\delta$ such that $\mathcal H\oplus\{0\}$ is invariant
under $X+\delta X^{-1}$ and
\[
 T=(X+\delta X^{-1})\big|_{\mathcal H\oplus\{0\}}.
\]
\end{theorem}

To apply the quantum annulus dilation of McCullough and Pascoe, put
$S=X/\sqrt\delta$. Then $\norm{S},\norm{S^{-1}}\leq1/\sqrt\delta$.
Their dilation $\widetilde S$ satisfies
$\widetilde S^*\widetilde S=\delta P+\delta^{-1}(I-P)$ and
$S^n=J^*\widetilde S^nJ$ for $n\in\Z$. Rescaling by
$U=\sqrt\delta\,\widetilde S$ and setting $E=I-P$ gives the following form.

\begin{theorem}[{\cite[Theorem~1.1(c)]{McCulloughPascoe}}]\label{thm:mp}
Let $0<\delta<1$, and let $X$ be an invertible operator on a Hilbert space
$\mathcal K_0$ with $\norm{X}\leq1$ and $\norm{\delta X^{-1}}\leq1$. Then
there are a Hilbert space $\mathcal K$, an isometry
$J:\mathcal K_0\to\mathcal K$, an orthogonal projection $E$ on
$\mathcal K$ and an invertible operator $U\in\mathcal B(\mathcal K)$ such
that
\begin{equation}\label{eq:laurent}
 \abs{U}=\delta I+(1-\delta)E,
 \qquad
 X^n=J^*U^nJ,
 \qquad
 n\in\Z.
\end{equation}
In particular, $\delta I\leq\abs{U}\leq I$, so that $\norm{U}\leq1$ and
$\norm{\delta U^{-1}}\leq1$.
\end{theorem}

Theorems~\ref{thm:aly} and~\ref{thm:mp} combine into a dilation of the
Faber calculus. For a polynomial symbol $g$,
Lemma~\ref{lem:faber-transform}{\rm(i)} shows that $f_g$ is a
polynomial, so $f_g(T)$ is defined for every operator $T$.

\begin{proposition}\label{prop:faber-dilation}
Let $0<\delta<1$, $\dim\mathcal H<\infty$ and $W(T)\subseteq K_\delta$.
Then there are a Hilbert space $\mathcal K$, an isometry
$J:\mathcal H\to\mathcal K$ and an invertible operator
$U\in\mathcal B(\mathcal K)$ with $\abs{U}=\delta I+(1-\delta)E$ for an
orthogonal projection $E$, such that
\begin{equation}\label{eq:faber-dilation}
 f_g(T)
 =
 J^*\bigl(g(U)+g(\delta U^{-1})-g(0)I\bigr)J
\end{equation}
for every polynomial $g$. Consequently,
\[
 \norm{f_g(T)+tI}
 \leq
 \norm{g(U)+g(\delta U^{-1})-g(0)I+tI},
 \qquad
 t\in\C.
\]
\end{proposition}

\begin{proof}
Let $X$ be given by Theorem~\ref{thm:aly}, let $Y=X+\delta X^{-1}$, and
write $\iota x=(x,0)$. Since $Y\iota=\iota T$, the same intertwining
identity holds for the polynomial $f_g$. Substitution of $X$ in the
Laurent polynomial identity~\eqref{eq:faber-poly} gives
\[
 F_{n,\delta}(Y)=X^n+\delta^nX^{-n}\quad(n\geq1),
 \qquad
 f_g(Y)=g(X)+g(\delta X^{-1})-g(0)I.
\]
Let $J_1$ and $U$ be the isometry and dilation of Theorem~\ref{thm:mp}
for $X$. Compressing each Laurent power by~\eqref{eq:laurent} gives
\[
 f_g(T)=\iota^*f_g(Y)\iota
       =(J_1\iota)^*\bigl(g(U)+g(\delta U^{-1})-g(0)I\bigr)(J_1\iota).
\]
This proves~\eqref{eq:faber-dilation} with $J=J_1\iota$. The norm
estimate follows by compressing the shifted operator with the same isometry.
\end{proof}

\subsection{The Faber calculus up to the boundary}

For a polynomial symbol, compression to a finite-dimensional subspace
preserves the value of $f_g(T)$ at any prescribed vector. To include the
disk, we approximate it by ellipses with $\delta\downarrow0$.

For $T\in\BH$, a unit vector $x$ and $N\geq1$, let
\[
 \mathcal H_N=\operatorname{span}\{x,Tx,\dots,T^Nx\},
 \qquad
 T_N=P_{\mathcal H_N}T\big|_{\mathcal H_N},
\]
where $P_{\mathcal H_N}$ is the orthogonal projection onto $\mathcal H_N$.
Then
\[
 \dim\mathcal H_N\leq N+1,
 \qquad
 W(T_N)\subseteq W(T).
\]

\begin{lemma}\label{lem:approx}
Let $g$ be a polynomial of degree at most $N$.

\noindent{\rm (i)} If $W(T)\subseteq K_\delta$ and $x$ is a unit vector,
then $W(T_N)\subseteq K_\delta$ and
\[
 f_g(T_N)x=f_g(T)x.
\]

\noindent{\rm (ii)} Let $W(T)\subseteq\overline\D$, and put
$T_\epsilon=(1-\epsilon)T$ for $0<\epsilon<1$. Then
$W(T_\epsilon)\subseteq K_\epsilon$, and
\[
 (\mathcal F_\epsilon g)(T_\epsilon)\longrightarrow g(T)
\]
in norm as $\epsilon\downarrow0$.
\end{lemma}

\begin{proof}
For part~{\rm (i)}, the inclusion $W(T_N)\subseteq K_\delta$ follows from
$W(T_N)\subseteq W(T)$. If $0\leq j<N$ and $T_N^jx=T^jx$, then
\[
 T_N^{j+1}x=P_{\mathcal H_N}T(T_N^jx)=P_{\mathcal H_N}T^{j+1}x=T^{j+1}x,
\]
because $T^{j+1}x\in\mathcal H_N$. By induction, $T_N^jx=T^jx$ for
$0\leq j\leq N$. Since $f_g$ is a polynomial of degree at most $N$ by
Lemma~\ref{lem:faber-transform}{\rm(i)}, the claim follows.

For part~{\rm (ii)}, note that
\[
 m_\epsilon(z)\leq\frac{\abs{z}}{1-\epsilon}\leq1,
 \qquad
 \abs{z}\leq1-\epsilon.
\]
Thus $W(T_\epsilon)\subseteq K_\epsilon$. By the recurrence for the Faber
polynomials, the coefficients of $F_{n,\epsilon}$ are polynomials in
$\epsilon$, and $F_{n,0}(z)=z^n$. Consequently, the coefficients of
$\mathcal F_\epsilon g$ converge to those of $g$ as $\epsilon\downarrow0$.
Since also $T_\epsilon\to T$ in norm, the claim follows.
\end{proof}

The following estimate gives a continuous calculus for disk-algebra
symbols even when the spectrum of $T$ meets $\partial K_\delta$.

\begin{proposition}\label{prop:boundary-calculus}
Let $W(T)\subseteq K_\delta$.

\noindent{\rm (i)} For every polynomial $h$,
\[
 \norm{f_h(T)}\leq3\norm{h}_\infty.
\]

\noindent{\rm (ii)} The map $h\mapsto f_h(T)$ extends uniquely to a
bounded linear map on $A(\overline\D)$, still denoted by
$h\mapsto f_h(T)$, with $\norm{f_h(T)}\leq3\norm{h}_\infty$. If
$\sigma(T)\subset G_\delta$, this extension agrees with the holomorphic
functional calculus of $f_h\in\Hol(G_\delta)$.
\end{proposition}

\begin{proof}
For a polynomial $h$, the estimate of Beckermann and Reichel
in~\eqref{eq:br} gives
\[
 \norm{f_h(T)}\leq \norm{f_h(T)+h(0)I}+\abs{h(0)}
 \leq 3\norm{h}_\infty.
\]
For arbitrary $\mathcal H$ the same estimate follows by applying the
matrix bound to $T_N$ in Lemma~\ref{lem:approx}{\rm(i)} at each unit vector.

For part~{\rm (ii)}, polynomials are dense in $A(\overline\D)$, so
part~{\rm (i)} gives a unique bounded extension. Suppose
$\sigma(T)\subset G_\delta$. Let $h\in A(\overline\D)$, and let $h_j$ be
polynomials with $h_j\to h$ uniformly on $\overline\D$. By
Lemma~\ref{lem:faber-transform}{\rm(ii)},
\[
 \sup_{G_\delta}\abs{f_{h_j}-f_h}\leq3\norm{h_j-h}_\infty\longrightarrow0.
\]
Thus $f_{h_j}\to f_h$ uniformly on a neighborhood of $\sigma(T)$, and
the continuity of the holomorphic calculus gives $f_{h_j}(T)\to f_h(T)$.
Therefore the two definitions of $f_h(T)$ agree.
\end{proof}

To pass from polynomial symbols to general symbols, we use the Fejér
means. For $g(u)=\sum_{j\geq0}a_ju^j$ in $H^\infty(\D)$ and $n\geq0$, put
\begin{equation}\label{eq:fejer}
 g_n(u)=\sum_{j=0}^n\Bigl(1-\frac{j}{n+1}\Bigr)a_ju^j.
\end{equation}
Then $g_n(0)=g(0)$. On each circle $\{\abs{u}=\varrho\}$ with
$0<\varrho<1$, the function $g_n$ is the convolution of
$\theta\mapsto g(\varrho e^{i\theta})$ with the Fejér kernel, which is
positive and has integral one; see~\cite[Chapter~I]{Katznelson}. Hence
\[
 \norm{g_n}_\infty\leq\norm{g}_\infty,
\]
and $g_n\in\Sch$ whenever $g\in\Sch$. If $g\in A(\overline\D)$, Fejér's
theorem gives $g_n\to g$ uniformly on $\T$, hence on $\overline\D$ by the
maximum principle. For every $g\in H^\infty(\D)$, one has $g_n\to g$
locally uniformly on $\D$, because the coefficients $a_j$ are bounded.

The two calculi lead to the following class of pairs.

\begin{definition}\label{def:admissible}
Let $0\leq s\leq1$. An \emph{admissible pair} is a pair $(g,T)$ with
$g\in\Sch$ and $T\in\BH$, for some Hilbert space $\mathcal H$, such that
$W(T)\subseteq K_\delta$ and either $g\in A(\overline\D)$ or
$\sigma(T)\subset G_\delta$. For $t\in\C$, we put
\[
 \Mfr(t)
 =
 \sup\bigl\{\norm{f_g(T)+tI}:(g,T)\ \text{admissible}\bigr\}.
\]
\end{definition}

In the first case of Definition~\ref{def:admissible}, $f_g(T)$ is defined
by Proposition~\ref{prop:boundary-calculus}{\rm(ii)}. In the second case,
it is defined by the holomorphic functional calculus. When both apply,
the two definitions agree by the same proposition.

The last result of this section shows that, to bound $\Mfr(t)$, it is
enough to consider matrices, $\delta>0$ and polynomial symbols.

\begin{lemma}\label{lem:reduction}
Let $0\leq s<1$ and $t\in\C$, and let $\mu:[0,1)\to[0,\infty)$ be
continuous at $0$. Assume that
\[
 \norm{f_g(T)+tI}\leq\mu(\delta)
\]
whenever $0<\delta<1$, $\dim\mathcal H<\infty$, $W(T)\subseteq K_\delta$
and $g\in\Sch$ is a polynomial. Then $\Mfr(t)\leq\mu(\delta)$ for every
$0\leq\delta<1$.
\end{lemma}

\begin{proof}
Fix $0\leq\delta<1$ and an admissible pair $(g,T)$ for $K_\delta$.

\smallskip
\noindent\emph{Step 1: polynomial $g$ and $\delta>0$.} Let $g$ have degree
at most $N$. For each unit vector $x$, Lemma~\ref{lem:approx}{\rm(i)} and
the hypothesis give
\[
 \norm{(f_g(T)+tI)x}
 =
 \norm{(f_g(T_N)+tI)x}
 \leq
 \mu(\delta).
\]

\smallskip
\noindent\emph{Step 2: polynomial $g$ and $\delta=0$.} Step~1, applied to
$T_\epsilon$ and $K_\epsilon$, gives
\[
 \norm{(\mathcal F_\epsilon g)(T_\epsilon)+tI}\leq\mu(\epsilon),
 \qquad
 0<\epsilon<1.
\]
By Lemma~\ref{lem:approx}{\rm(ii)} and the continuity of $\mu$ at $0$,
\[
 \norm{g(T)+tI}\leq\lim_{\epsilon\downarrow0}\mu(\epsilon)=\mu(0).
\]

\smallskip
\noindent\emph{Step 3: general $g$.} The Fejér means~\eqref{eq:fejer} are
polynomials in $\Sch$. Hence Steps~1 and~2 give
\[
 \norm{f_{g_n}(T)+tI}\leq\mu(\delta),
 \qquad
 n\geq0.
\]
If $g\in A(\overline\D)$, then $g_n\to g$ uniformly on $\overline\D$, and
Proposition~\ref{prop:boundary-calculus}{\rm(ii)} gives
$f_{g_n}(T)\to f_g(T)$.

If $\sigma(T)\subset G_\delta$, then $g_n\to g$ locally uniformly on
$\D$. By Lemma~\ref{lem:levels}, the roots $u_1,u_2$ stay in a compact
subset of $\D$ when $z$ ranges over a compact subset of $G_\delta$. Therefore
$f_{g_n}\to f_g$ locally uniformly on $G_\delta$, and the continuity of
the holomorphic calculus gives $f_{g_n}(T)\to f_g(T)$.

In both cases, the bound passes to the limit.
\end{proof}

\section{The sharp shifted norm}\label{sec:block}

In this section we prove Theorem~\ref{thm:main}. Put
\begin{equation}\label{eq:M}
 M(t)
 =
 M_{\delta,s}(t)
 =
 \max\Bigl\{R+\abs{C+t},\ r+\sqrt{d^2+\abs{c+t}^2}\Bigr\},
 \qquad
 t\in\C.
\end{equation}
By Proposition~\ref{prop:faber-dilation} and Lemma~\ref{lem:reduction},
it suffices to estimate the operator $g(U)+g(\delta U^{-1})-sI$, where
$\abs{U}$ takes only the values $\delta$ and $1$.

In Subsections~\ref{ssec:one-sided}--\ref{ssec:upper}, we assume
$0<\delta<1$, $0\leq s<1$ and $k=s^2$. Let $U$ be
an invertible operator on a Hilbert space $\mathcal K$ with
\begin{equation}\label{eq:standing}
 \abs{U}=\delta I+(1-\delta)E_1
\end{equation}
for an orthogonal projection $E_1$, as in Theorem~\ref{thm:mp}. We write
\[
 B=\abs{U},
 \qquad
 D=\abs{U^*},
 \qquad
 V=\delta U^{-1}.
\]
Let $U=\Theta B$ be the polar decomposition; $\Theta$ is unitary because
$U$ is invertible. Then
\[
 U=D\Theta,
 \qquad
 D=\Theta B\Theta^*=\delta I+(1-\delta)E_2,
 \qquad
 E_2=\Theta E_1\Theta^*.
\]
Put $B^\sharp=\delta B^{-1}$ and $D^\sharp=\delta D^{-1}$. Since the spectra
of $B$ and $D$ lie in $\{\delta,1\}$,
\begin{equation}\label{eq:sharp}
 B^\sharp=(1+\delta)I-B,
 \qquad
 D^\sharp=(1+\delta)I-D.
\end{equation}
Moreover, $V=\delta B^{-1}\Theta^*=\Theta^*\delta D^{-1}$, that is,
\[
 V=B^\sharp\Theta^*=\Theta^*D^\sharp.
\]
Thus $\abs{V}=D^\sharp$, $\abs{V^*}=B^\sharp$, and
$\norm{U},\norm{V}\leq1$.

The constants~\eqref{eq:parameters} are endpoint values of three
functions on the interval $[2\delta,1+\delta^2]$:
\begin{equation}\label{eq:functions}
 e(v)=1-kv+k^2\delta^2,
 \qquad
 c(v)=s-s(1-k)\frac{v-2k\delta^2}{e(v)},
 \qquad
 r(v)=\frac{(1-k)(1+\delta)(1-k\delta)}{e(v)}.
\end{equation}
Indeed, $e(2\delta)=(1-k\delta)^2$ and $e(1+\delta^2)=(1-k)(1-k\delta^2)$,
and a direct substitution gives
\begin{equation}\label{eq:endpoints}
 c(2\delta)=c,
 \quad
 r(2\delta)=r,
 \quad
 \frac{1-\delta}{1+\delta}\,r(2\delta)=d,
 \qquad
 c(1+\delta^2)=C,
 \quad
 r(1+\delta^2)=R.
\end{equation}
The function $e$ is positive on $[2\delta,1+\delta^2]$, since it is affine
and positive at both endpoints.

\subsection{A one-sided inequality}\label{ssec:one-sided}

We use the Heinz inequality~\cite{Heinz}: for positive operators $P,Q$
and any operator $Y$,
\begin{equation}\label{eq:heinz}
 \norm{P^{1/2}YQ^{1/2}}\leq\norm{PYQ}^{1/2}\norm{Y}^{1/2};
\end{equation}
see~\cite{FujiiFujiiFurutaNakamoto} for an elementary proof.

\begin{lemma}\label{lem:one-sided}
Let $X$ be an invertible operator with $\norm{X}\leq1$, and write
$X=\Theta_0B_0=D_0\Theta_0$ with $\Theta_0$ unitary, $B_0=\abs{X}$ and
$D_0=\abs{X^*}$; thus $B_0,D_0$ are positive, invertible and at most $I$.
Let $g\in\Sch\cap\Hol(\overline\D)$, and define $Y$ by
\[
 g(X)=sI+(1-k)Y.
\]
With
\[
 S_D=(D_0^{-1}-kB_0)^{-1},
 \qquad
 S_B=(B_0^{-1}-kD_0)^{-1},
\]
one has
\begin{equation}\label{eq:one-sided}
 \begin{pmatrix}S_D&Y+sS_DB_0\\(Y+sS_DB_0)^*&S_B\end{pmatrix}\geq0.
\end{equation}
\end{lemma}

\begin{proof}
\emph{Step 1: a Schwarz inequality.} The function
$\widetilde g=(g-s)/(1-sg)$ is holomorphic on a neighborhood of
$\overline\D$, maps $\D$ into $\overline\D$, and vanishes at $0$. By the
Schwarz lemma,
\[
 \widetilde g(u)=u\psi(u),
\]
where $\psi$ is holomorphic on a neighborhood of $\overline\D$ and bounded
by one on $\D$. Put
\[
 G=\widetilde g(X)=X\psi(X)=\psi(X)X.
\]
Von Neumann's inequality~\cite{vonNeumann} gives $\norm{\psi(X)}\leq1$. Since
\[
 D_0^{-1}G=\Theta_0\psi(X),
 \qquad
 GB_0^{-1}=\psi(X)\Theta_0,
\]
both operators are contractions. We apply~\eqref{eq:heinz} with
\[
 P=D_0^{-1},
 \qquad
 Q=B_0,
 \qquad
 Z_0=GB_0^{-1}
\]
in place of $Y$. Since $PZ_0Q=D_0^{-1}G$, and both $PZ_0Q$ and $Z_0$ are
contractions, we obtain
\[
 \norm{D_0^{-1/2}GB_0^{-1/2}}\leq1.
\]
Conjugation by $B_0^{1/2}$ gives
\begin{equation}\label{eq:schwarz}
 G^*D_0^{-1}G\leq B_0.
\end{equation}

\smallskip
\noindent\emph{Step 2: the Möbius change of variable.} Since
$\norm{G}\leq1$ and $s<1$, the operator $I+sG$ is invertible. From
$g=(s+\widetilde g)/(1+s\widetilde g)$, we get
\[
 g(X)=(sI+G)(I+sG)^{-1},
\]
and therefore
\[
 Y=G(I+sG)^{-1},
 \qquad
 I-sY=(I+sG)^{-1},
 \qquad
 G(I-sY)=Y.
\]
Conjugating~\eqref{eq:schwarz} by $I-sY$ gives
$Y^*D_0^{-1}Y\leq(I-sY)^*B_0(I-sY)$, that is,
\begin{equation}\label{eq:expanded}
 Y^*(D_0^{-1}-kB_0)Y+s(Y^*B_0+B_0Y)\leq B_0.
\end{equation}

\smallskip
\noindent\emph{Step 3: completing the square.} Since $D_0\leq I$ and
$B_0\leq I$,
\[
 D_0^{-1}-kB_0\geq(1-k)I,
 \qquad
 B_0^{-1}-kD_0\geq(1-k)I.
\]
Hence $S_D$ and $S_B$ are positive and invertible. Since
$S_D^{-1}=D_0^{-1}-kB_0$, expanding the square and
using~\eqref{eq:expanded}, we obtain
\begin{align*}
 (Y+sS_DB_0)^*S_D^{-1}(Y+sS_DB_0)
 &=
 Y^*S_D^{-1}Y+s(Y^*B_0+B_0Y)+kB_0S_DB_0\\
 &\leq
 B_0+kB_0S_DB_0.
\end{align*}

The right-hand side equals $S_B$. Indeed,
$(I-kD_0B_0)S_D=D_0(D_0^{-1}-kB_0)S_D=D_0$, and hence
\[
 (B_0^{-1}-kD_0)(B_0+kB_0S_DB_0)
 =
 I-kD_0B_0+k(I-kD_0B_0)S_DB_0
 =
 I.
\]
Therefore
\[
 (Y+sS_DB_0)^*S_D^{-1}(Y+sS_DB_0)\leq S_B.
\]
By the Schur complement criterion for positive block
operators~\cite[Theorem~1.3.3]{Bhatia}, this is equivalent
to~\eqref{eq:one-sided}.
\end{proof}

\subsection{The two-sided inequality}

We apply Lemma~\ref{lem:one-sided} to $X=U$, with
$(\Theta_0,B_0,D_0)=(\Theta,B,D)$, and to $X=V$, with
$(\Theta_0,B_0,D_0)=(\Theta^*,D^\sharp,B^\sharp)$, and we add the
results. Put
\[
 \Omega_1=DB,
 \qquad
 \Omega_2=B^\sharp D^\sharp,
 \qquad
 H=\Omega_1+\Omega_2,
 \qquad
 \Sigma=\frac{D-B}{1+\delta},
 \qquad
 Z=g(U)+g(V)-sI.
\]

Adding the inequalities for $U$ and $V$ expresses both diagonal bounds
and the center of the off-diagonal entry as functions of $\Sigma$.
This will allow us to use its scalar spectral variable.

\begin{proposition}\label{prop:block}
Let $g\in\Sch\cap\Hol(\overline\D)$.

\noindent{\rm (i)} The operator $\Sigma$ is self-adjoint, with
$\norm{\Sigma}\leq j_0=(1-\delta)/(1+\delta)$, and
\[
 H=(1+\delta^2)I-(1+\delta)^2\Sigma^2,
 \qquad
 \Omega_1\Omega_2=\Omega_2\Omega_1=\delta^2I.
\]
Moreover, $H$ commutes with $B$, $D$, $B^\sharp$, $D^\sharp$, $\Omega_1$
and $\Omega_2$.

\noindent{\rm (ii)} With the functions~\eqref{eq:functions},
\begin{equation}\label{eq:block}
 \begin{pmatrix}r(H)(I+\Sigma)&Z-c(H)\\(Z-c(H))^*&r(H)(I-\Sigma)\end{pmatrix}
 \geq0.
\end{equation}
\end{proposition}

\begin{proof}
For part~{\rm (i)}, write
\[
 B=\delta I+(1-\delta)E_1,
 \qquad
 D=\delta I+(1-\delta)E_2.
\]
Then $\Sigma=j_0(E_2-E_1)$, and $E_2-E_1$ is self-adjoint with norm at
most one. By~\eqref{eq:sharp},
\[
 B^\sharp=I-(1-\delta)E_1,
 \qquad
 D^\sharp=I-(1-\delta)E_2.
\]
Expanding, we obtain
\begin{align*}
 H
 &=
 (1+\delta^2)I-(1-\delta)^2(E_1+E_2-E_1E_2-E_2E_1)\\
 &=
 (1+\delta^2)I-(1-\delta)^2(E_2-E_1)^2,
\end{align*}
which is the stated formula. Also,
\[
 \Omega_1\Omega_2=DBB^\sharp D^\sharp=\delta DD^\sharp=\delta^2I,
\]
and similarly $\Omega_2\Omega_1=\delta^2I$.

A direct computation gives
\[
 E_1(E_2-E_1)^2=E_1-E_1E_2E_1=(E_2-E_1)^2E_1,
\]
and similarly for $E_2$. Thus $H$ commutes with $E_1$ and $E_2$, and
therefore with $B$, $D$, $B^\sharp$, $D^\sharp$, $\Omega_1$ and $\Omega_2$.

For part~{\rm (ii)}, define $Y_1$ and $Y_2$ by
\[
 g(U)=sI+(1-k)Y_1,
 \qquad
 g(V)=sI+(1-k)Y_2,
\]
so that
\[
 Z-sI=(1-k)(Y_1+Y_2).
\]
Lemma~\ref{lem:one-sided} gives two positive block operators. Their upper diagonal entries are
\[
 (D^{-1}-kB)^{-1}=(I-k\Omega_1)^{-1}D,
 \qquad
 \bigl((B^\sharp)^{-1}-kD^\sharp\bigr)^{-1}=(I-k\Omega_2)^{-1}B^\sharp.
\]
Since $\Omega_1$ and $\Omega_2$ commute, part~{\rm (i)} gives
\begin{align*}
 (I-k\Omega_1)(I-k\Omega_2)
 &=
 I-kH+k^2\delta^2I
 =
 e(H)\\
 &\geq
 (1-k)(1-k\delta^2)I
 >0.
\end{align*}
Moreover, $e(H)$ commutes with $B$, $D$, $B^\sharp$ and $D^\sharp$.

Using $\Omega_2D=\delta B^\sharp$ and $\Omega_1B^\sharp=\delta D$, we rewrite
the upper diagonal entries as
\[
 (I-k\Omega_1)^{-1}D=(D-k\delta B^\sharp)e(H)^{-1},
 \qquad
 (I-k\Omega_2)^{-1}B^\sharp=(B^\sharp-k\delta D)e(H)^{-1}.
\]
For the lower diagonal entries, we use that $H$ is self-adjoint, so that
also $H=BD+D^\sharp B^\sharp$. Since moreover
$BD\,D^\sharp B^\sharp=\delta^2I$, one has
$(I-kBD)(I-kD^\sharp B^\sharp)=e(H)$, and therefore
\begin{align*}
 (B^{-1}-kD)^{-1}
 &=
 (I-kBD)^{-1}B
 =
 (I-kD^\sharp B^\sharp)e(H)^{-1}B\\
 &=
 (B-k\delta D^\sharp)e(H)^{-1}.
\end{align*}
In the same way,
\[
 \bigl((D^\sharp)^{-1}-kB^\sharp\bigr)^{-1}=(D^\sharp-k\delta B)e(H)^{-1}.
\]
The off-diagonal entries are $Y_1+s(I-k\Omega_1)^{-1}\Omega_1$ and
$Y_2+s(I-k\Omega_2)^{-1}\Omega_2$, and
\begin{align*}
 (I-k\Omega_1)^{-1}\Omega_1+(I-k\Omega_2)^{-1}\Omega_2
 &=
 \bigl((I-k\Omega_2)\Omega_1+(I-k\Omega_1)\Omega_2\bigr)e(H)^{-1}\\
 &=
 (H-2k\delta^2I)e(H)^{-1}.
\end{align*}

We now add the two blocks and multiply by $1-k$. By~\eqref{eq:sharp},
\[
 D+B^\sharp=(1+\delta)(I+\Sigma),
 \qquad
 B+D^\sharp=(1+\delta)(I-\Sigma).
\]
Hence the diagonal entries become
\[
 (1-k)(1-k\delta)(D+B^\sharp)e(H)^{-1}=r(H)(I+\Sigma),
 \qquad
 (1-k)(1-k\delta)(B+D^\sharp)e(H)^{-1}=r(H)(I-\Sigma),
\]
and the off-diagonal entry becomes
\[
 Z-sI+s(1-k)(H-2k\delta^2I)e(H)^{-1}=Z-c(H).
\]
This proves~\eqref{eq:block}.
\end{proof}

\subsection{The scalar inequality and the upper bound}\label{ssec:upper}

By Proposition~\ref{prop:block}{\rm(i)}, the
operator $H$ is the function
\[
 v(j)=1+\delta^2-(1+\delta)^2j^2,
 \qquad
 -j_0\leq j\leq j_0,
\]
of $\Sigma$, and $v$ maps $[-j_0,j_0]$ onto $[2\delta,1+\delta^2]$. For
$t\in\C$ and $-j_0\leq j\leq j_0$, put
\[
 \Xi_t(j)
 =
 \begin{pmatrix}r(v)(1+j)&c(v)+t\\c(v)+\overline t&r(v)(1-j)\end{pmatrix},
 \qquad
 v=v(j).
\]

\begin{lemma}\label{lem:scalar}
For every $t\in\C$ and $-j_0\leq j\leq j_0$,
\[
 \Xi_t(j)\leq M(t)I_2.
\]
\end{lemma}

\begin{proof}
Put $M=M(t)$ and $v=v(j)$. A Hermitian $2\times2$ matrix is positive if
and only if its trace and its determinant are nonnegative. The trace of
$MI_2-\Xi_t(j)$ is $2(M-r(v))$, and its determinant is
\[
 Q(v)
 =
 \bigl(M-r(v)\bigr)^2-r(v)^2j^2-\abs{c(v)+t}^2.
\]

We first check the trace. The function $r$ is a positive constant
divided by the positive affine function $e$. Thus its maximum on
$[2\delta,1+\delta^2]$ is attained at an endpoint, and~\eqref{eq:endpoints}
gives
\[
 r(v)\leq\max\{r,R\}\leq M.
\]

For the determinant, we use $1-j^2=(v+2\delta)/(1+\delta)^2$ and expand:
\[
 Q(v)
 =
 M^2-\abs{t}^2-2Mr(v)-2c(v)\Real t
 +r(v)^2\frac{v+2\delta}{(1+\delta)^2}-c(v)^2.
\]
The last two terms combine into a single fraction with denominator
$e(v)$. Put
\[
 x=v+2\delta,
 \qquad
 \gamma_1=1+k\delta,
 \qquad
 \gamma_2=1+(2-k)\delta.
\]
Then
\[
 e(v)=\gamma_1^2-kx,
 \qquad
 (1-k)(1-k\delta)=\gamma_1-k\gamma_2,
 \qquad
 c(v)=\frac{s(\gamma_1\gamma_2-x)}{e(v)}.
\]
The identity
\[
 (\gamma_1-k\gamma_2)^2x-k(\gamma_1\gamma_2-x)^2
 =
 (\gamma_1^2-kx)(x-k\gamma_2^2)
\]
holds because both sides are quadratic in $x$, with leading coefficient
$-k$ and constant term $-k\gamma_1^2\gamma_2^2$, and their coefficients of
$x$ agree:
\[
 (\gamma_1-k\gamma_2)^2+2k\gamma_1\gamma_2=\gamma_1^2+k^2\gamma_2^2.
\]
Since $k=s^2$, this identity gives
\begin{equation}\label{eq:key}
 r(v)^2\frac{v+2\delta}{(1+\delta)^2}-c(v)^2
 =
 \frac{v+2\delta-k\gamma_2^2}{e(v)}.
\end{equation}
Consequently, every term of $Q(v)$ is either a constant or an affine
function of $v$ divided by $e(v)$. Hence $e(v)Q(v)$ is an affine function
of $v$.

At the endpoints,~\eqref{eq:endpoints} gives
\[
 Q(2\delta)=(M-r)^2-d^2-\abs{c+t}^2,
 \qquad
 Q(1+\delta^2)=(M-R)^2-\abs{C+t}^2.
\]
Both are nonnegative by the definition of $M$. Since $e>0$, the affine
function $eQ$ is nonnegative at both endpoints, hence on the whole
interval $[2\delta,1+\delta^2]$. Therefore $Q(v)\geq0$, and
$MI_2-\Xi_t(j)\geq0$.
\end{proof}

\begin{theorem}\label{thm:upper-annulus}
Let $U$ satisfy~\eqref{eq:standing}, let $g\in\Sch\cap\Hol(\overline\D)$,
and let $t\in\C$. Then
\[
 \norm{g(U)+g(\delta U^{-1})-sI+tI}\leq M(t).
\]
\end{theorem}

\begin{proof}
Put $F=Z+tI$, where $Z=g(U)+g(V)-sI$, and consider the self-adjoint
operator
\[
 \widehat F=\begin{pmatrix}0&F\\F^*&0\end{pmatrix}
\]
on $\mathcal K\oplus\mathcal K$. One has $\norm{\widehat F}=\norm{F}$.
Conjugation by the unitary $\operatorname{diag}(I,-I)$ changes
$\widehat F$ into $-\widehat F$, so the spectrum of $\widehat F$ is
symmetric with respect to $0$. Thus $\norm{F}$ is the largest point of
$\sigma(\widehat F)$, and it suffices to show that
$\widehat F\leq M(t)I$.

Conjugating~\eqref{eq:block} by $\operatorname{diag}(I,-I)$, we obtain
\[
 \begin{pmatrix}0&Z-c(H)\\(Z-c(H))^*&0\end{pmatrix}
 \leq
 \begin{pmatrix}r(H)(I+\Sigma)&0\\0&r(H)(I-\Sigma)\end{pmatrix}.
\]
Since $c(H)$ is self-adjoint, we may add the self-adjoint operator
\[
 \begin{pmatrix}0&c(H)+tI\\c(H)+\overline tI&0\end{pmatrix}
\]
to both sides. This gives
\[
 \widehat F
 \leq
 \begin{pmatrix}r(H)(I+\Sigma)&c(H)+tI\\c(H)+\overline tI&r(H)(I-\Sigma)\end{pmatrix}
 =
 \Xi_t(\Sigma),
\]
where $\Xi_t(\Sigma)$ is obtained by applying the continuous functional
calculus of $\Sigma$ entrywise.

By Proposition~\ref{prop:block}{\rm(i)}, $\sigma(\Sigma)\subseteq[-j_0,j_0]$,
and Lemma~\ref{lem:scalar} gives $M(t)I_2-\Xi_t(j)\geq0$ for every
$j\in\sigma(\Sigma)$. The entrywise functional calculus is a
$*$-homomorphism from the $2\times2$ matrices over $C(\sigma(\Sigma))$ into
the $2\times2$ operator matrices on $\mathcal K$, and $*$-homomorphisms
preserve positivity. Hence $\Xi_t(\Sigma)\leq M(t)I$, and therefore
$\widehat F\leq M(t)I$.
\end{proof}

\begin{corollary}\label{cor:upper}
For $0\leq\delta<1$, $0\leq s<1$ and $t\in\C$,
\[
 \Mfr(t)\leq M_{\delta,s}(t).
\]
\end{corollary}

\begin{proof}
Let $0<\delta<1$, $\dim\mathcal H<\infty$ and $W(T)\subseteq K_\delta$, and
let $g\in\Sch$ be a polynomial. Proposition~\ref{prop:faber-dilation} and
Theorem~\ref{thm:upper-annulus} give
\[
 \norm{f_g(T)+tI}\leq M_{\delta,s}(t).
\]
The constants~\eqref{eq:parameters} are continuous in $\delta\in[0,1)$,
and so is $M_{\delta,s}(t)$. Hence Lemma~\ref{lem:reduction} with $\mu(\delta)=M_{\delta,s}(t)$ completes the
proof.
\end{proof}

\subsection{The extremal symbols and the two models}\label{ssec:models}

For each shift $t$, a disk automorphism suffices to attain each value
in~\eqref{eq:M}. In the rest of this section, $0\leq\delta<1$,
$0\leq s<1$ and $k=s^2$. The disk
automorphisms
\begin{equation}\label{eq:automorphism}
 \beta_s(u)=\frac{s+u}{1+su},
 \qquad
 g_\eta(u)=\beta_s(\eta u),
 \qquad
 \eta\in\T,
\end{equation}
belong to $\Sch\cap A(\overline\D)$, and their Faber transforms are
linear-fractional.

\begin{example}\label{ex:mobius}
Let $0\leq s<1$ and $\eta\in\T$. Since
\[
 g_\eta(u)-s=\frac{(1-k)\eta u}{1+s\eta u},
\]
the relations~\eqref{eq:roots} give
\begin{equation}\label{eq:mob-lift}
 f_{g_\eta}(z)
 =
 s+(1-k)\,\frac{\eta z+2s\delta\eta^2}{1+s\eta z+k\delta\eta^2}.
\end{equation}
The denominator equals
\[
 (1+s\eta u_1)(1+s\eta u_2),
\]
which does not vanish on the neighborhood $\{\varrho_\delta<1/s\}$ of
$K_\delta$, the whole plane if $s=0$.
\end{example}

To realize the two values in~\eqref{eq:M}, the phase $\eta$ must be
matched to the shift. For $\omega\in\T$, put
\begin{equation}\label{eq:phases}
 P_\omega=\frac{\omega-s}{1-s\omega},
 \qquad
 Q_\omega=\frac{\omega-s\delta}{1-s\delta\omega}.
\end{equation}
Both are unimodular, and $\beta_s(P_\omega)=\omega$.

\begin{lemma}[Scalar model]\label{lem:scalar-model}
Let $\omega\in\T$, choose $\eta\in\T$ with $\eta^2=P_\omega Q_\omega$, and
put
\[
 u=\frac{P_\omega}{\eta},
 \qquad
 \lambda=u+\delta\overline u.
\]
Then $\lambda\in\partial K_\delta$ and
\[
 f_{g_\eta}(\lambda)=C+R\omega.
\]
\end{lemma}

\begin{proof}
Since $\abs{u}=1$, one has $\lambda=L_\delta(u)\in\partial K_\delta$, and
the roots~\eqref{eq:roots} at $\lambda$ are $u$ and $\delta\overline u$.
Moreover,
\[
 \eta u=P_\omega,
 \qquad
 \eta\overline u=\frac{\eta^2}{P_\omega}=Q_\omega.
\]
Hence
\[
 f_{g_\eta}(\lambda)
 =
 \beta_s(\eta u)+\beta_s(\delta\eta\overline u)-s
 =
 \beta_s(P_\omega)+\beta_s(\delta Q_\omega)-s.
\]

The first term equals $\omega$. Clearing the denominator of $Q_\omega$, we
obtain
\begin{align*}
 \beta_s(\delta Q_\omega)
 &=
 \frac{s(1-\delta^2)+(1-k)\delta\omega}{1-k\delta^2}\\
 &=
 s-s\delta b+b\omega,
\end{align*}
with $b=(1-k)\delta/(1-k\delta^2)$ as in~\eqref{eq:parameters}. Therefore
\[
 f_{g_\eta}(\lambda)=(1+b)\omega-s\delta b=C+R\omega.
 \qedhere
\]
\end{proof}

\begin{lemma}[Focal model]\label{lem:focal-model}
Let $\omega\in\T$, and let $\mu,\nu\geq0$ with $\mu^2+\nu^2=1$. Choose
$v_1,v_2\in\T$ with
\begin{equation}\label{eq:phase}
 (\nu-s\delta\omega\mu)v_1=\omega\nu-s\delta\mu,
 \qquad
 (\mu-s\omega\nu)v_2=\omega\mu-s\nu,
\end{equation}
and $\eta\in\T$ with $\eta^2=v_1v_2$. Put
\[
 A_v=\begin{pmatrix}0&v_1\\\delta v_2&0\end{pmatrix},
 \qquad
 T_v=\frac{2}{\eta}A_v,
 \qquad
 x=\begin{pmatrix}\mu\\\nu\end{pmatrix},
 \qquad
 y=\begin{pmatrix}\nu\\\mu\end{pmatrix}.
\]
Then $T_v$ is unitarily equivalent to $A_\delta$, and
\begin{equation}\label{eq:focal-identity}
 \bigl(f_{g_\eta}(T_v)-cI\bigr)x
 =
 \omega\begin{pmatrix}(r+d)\nu\\(r-d)\mu\end{pmatrix}.
\end{equation}
If $\mu=\nu=1/\sqrt2$, one may take $v_1=Q_\omega$ and $v_2=P_\omega$. Then
$\eta^2=P_\omega Q_\omega$ as in Lemma~\ref{lem:scalar-model}, and
\[
 \inner{f_{g_\eta}(T_v)x}{x}=c+r\omega.
\]
If moreover $\omega=1$, one may take $\eta=1$. Then $T_v=A_\delta$ and
\[
 f_{g_1}(A_\delta)=cI+aA_\delta,
\]
with $a=(1-k)/(1-k\delta)$ as in~\eqref{eq:parameters}.
\end{lemma}

\begin{proof}
In each equation~\eqref{eq:phase}, the two coefficients have equal
modulus, because $\abs{\omega}=1$ and $s,\delta,\mu,\nu$ are real. If they
are nonzero, $v_1$ and $v_2$ are determined; otherwise any unimodular
value works.

Put $\theta=v_1/\eta$. Then $\theta\in\T$ and $v_2/\eta=\eta/v_1=\overline\theta$.
Hence
\[
 T_v
 =
 \begin{pmatrix}0&2\theta\\2\delta\overline\theta&0\end{pmatrix}
 =
 \Pi A_\delta\Pi^*,
 \qquad
 \Pi=\operatorname{diag}(\theta,1),
\]
so $T_v$ is unitarily equivalent to $A_\delta$.

The matrix $U_0=A_v/\eta$ satisfies $U_0^2=\delta I$ and $T_v=2U_0$. For
$\delta>0$, one has $T_v=U_0+\delta U_0^{-1}$, and~\eqref{eq:faber-poly}
gives
\[
 F_{n,\delta}(T_v)=U_0^n+\delta^nU_0^{-n}=2U_0^n,
 \qquad
 n\geq1;
\]
for $\delta=0$, the same identity holds because $U_0^2=0$. Since
$\sigma(T_v)\subset G_\delta$, Lemma~\ref{lem:faber-transform}{\rm(i)}
allows us to sum the Faber series of $g_\eta$ at $T_v$. As
$g_\eta(U_0)=\beta_s(A_v)$, we obtain
\[
 Z_0:=f_{g_\eta}(T_v)=2\beta_s(A_v)-sI=sI+2(1-k)A_v(I+sA_v)^{-1}.
\]

We now subtract $cI$. Since $s-c=2s\delta a$ and $1-k=a(1-k\delta)$,
\begin{align*}
 Z_0-cI
 &=
 (I+sA_v)^{-1}\bigl(2(1-k)A_v+2s\delta a(I+sA_v)\bigr)\\
 &=
 2a(I+sA_v)^{-1}(A_v+s\delta I).
\end{align*}
Comparing coordinates, the two equations~\eqref{eq:phase} state exactly
that
\[
 (A_v+s\delta I)x=\omega(I+sA_v)\Delta y,
 \qquad
 \Delta=\operatorname{diag}(1,\delta).
\]
Since $2a\Delta=\operatorname{diag}(r+d,r-d)$ by~\eqref{eq:parameters},
this gives~\eqref{eq:focal-identity}.

Finally, let $\mu=\nu$. Then~\eqref{eq:phase} reads
\[
 (1-s\delta\omega)v_1=\omega-s\delta,
 \qquad
 (1-s\omega)v_2=\omega-s,
\]
so that $v_1=Q_\omega$ and $v_2=P_\omega$, and~\eqref{eq:focal-identity}
gives
\[
 \inner{(Z_0-cI)x}{x}=\frac{\omega}{2}\bigl((r+d)+(r-d)\bigr)=r\omega.
\]
If moreover $\omega=1$, then $v_1=v_2=1$, and we may take $\eta=1$. In
this case $T_v=A_\delta$ and $A_v=A_\delta/2$. Since $A_v^2=\delta I$, one
has $(I+sA_v)A_v=A_v+s\delta I$, and hence
\[
 Z_0-cI=2aA_v=aA_\delta.
 \qedhere
\]
\end{proof}

\subsection{Attainment and the common phase}\label{ssec:attainment}

The models of Lemmas~\ref{lem:scalar-model} and~\ref{lem:focal-model}
give the following attainment statement for Theorem~\ref{thm:main}.

\begin{theorem}\label{thm:shift}
Let $0\leq\delta<1$ and $0\leq s<1$. For every $t\in\C$,
\begin{equation}\label{eq:main}
 \Mfr(t)=M_{\delta,s}(t),
\end{equation}
with $M_{\delta,s}$ as in~\eqref{eq:M}. The first value in~\eqref{eq:M} is
attained by a disk automorphism $g_\eta$ and a scalar
$\lambda\in\partial K_\delta$, and the second by a disk automorphism
$g_\eta$ and a matrix unitarily equivalent to $A_\delta$. In each model,
the shifted norm equals the corresponding value in~\eqref{eq:M}, also
when the other value is larger.
\end{theorem}

\begin{proof}
Corollary~\ref{cor:upper} gives $\Mfr(t)\leq M_{\delta,s}(t)$. It remains to
attain the two values. All pairs below are admissible, because
$g_\eta\in A(\overline\D)$.

For the first value, let $\omega=(C+t)/\abs{C+t}$, or any $\omega\in\T$ if
$C+t=0$, and let $\eta$ and $\lambda$ be as in
Lemma~\ref{lem:scalar-model}. Since $C+t=\abs{C+t}\omega$,
\[
 \abs{f_{g_\eta}(\lambda)+t}
 =
 \abs{C+t+R\omega}
 =
 R+\abs{C+t}.
\]

For the second value, put
\begin{equation}\label{eq:mu-nu}
 \ell=\abs{c+t},
 \qquad
 \kappa=\sqrt{d^2+\ell^2},
 \qquad
 \mu=\sqrt{\frac{\kappa-d}{2\kappa}},
 \qquad
 \nu=\sqrt{\frac{\kappa+d}{2\kappa}},
\end{equation}
and let $\omega=(c+t)/\ell$, or any $\omega\in\T$ if $\ell=0$. Then
$\mu^2+\nu^2=1$, and
\[
 \ell\mu=(\kappa-d)\nu,
 \qquad
 \ell\nu=(\kappa+d)\mu.
\]
Let $\eta$, $T_v$, $x$ and $y$ be as in Lemma~\ref{lem:focal-model}. Since
$(c+t)x=\omega\ell x$, adding this to~\eqref{eq:focal-identity} gives
\[
 \bigl(f_{g_\eta}(T_v)+tI\bigr)x
 =
 \omega\begin{pmatrix}(r+d)\nu+\ell\mu\\(r-d)\mu+\ell\nu\end{pmatrix}
 =
 \omega(r+\kappa)\,y.
\]
Since $x$ and $y$ are unit vectors,
$\norm{f_{g_\eta}(T_v)+tI}\geq r+\kappa$.

For the reverse inequality when $\delta>0$, take $U_0=A_v/\eta$ as in
Lemma~\ref{lem:focal-model}. Since $\delta U_0^{-1}=U_0$,
\[
 f_{g_\eta}(T_v)=g_\eta(U_0)+g_\eta(\delta U_0^{-1})-sI.
\]
Here $\abs{U_0}=\operatorname{diag}(\delta,1)$ and
$\abs{U_0^*}=\operatorname{diag}(1,\delta)$, so that
$\Sigma=j_0\operatorname{diag}(1,-1)$. The spectral bound in the proof
of Theorem~\ref{thm:upper-annulus} therefore involves only
\[
 \Xi_t(\pm j_0)
 =
 \begin{pmatrix}r\pm d&c+t\\c+\overline t&r\mp d\end{pmatrix}.
\]
Each matrix has larger eigenvalue $r+\kappa$. Hence
$\norm{f_{g_\eta}(T_v)+tI}\leq r+\kappa$, even when the scalar value
is larger.

Suppose next that $\delta=0$. Then $A_v^2=0$, so that
$(I+sA_v)^{-1}=I-sA_v$ and $A_v(I+sA_v)^{-1}=A_v$. Moreover, $c=s$ and
$r=d=1-k$ by~\eqref{eq:parameters}. Hence
\[
 N:=f_{g_\eta}(T_v)+tI=(c+t)I+2dA_v,
 \qquad
 A_v=\begin{pmatrix}0&v_1\\0&0\end{pmatrix}.
\]
A direct computation gives
\[
 \operatorname{tr}(N^*N)=2\abs{c+t}^2+4d^2,
 \qquad
 \det(N^*N)=\abs{c+t}^4.
\]
Hence the largest eigenvalue of $N^*N$ is
\[
 \abs{c+t}^2+2d^2+2d\kappa=(d+\kappa)^2,
\]
and $\norm{N}=d+\kappa=r+\kappa$.

This completes the proof of~\eqref{eq:main}.
\end{proof}

At $t=s$, Theorem~\ref{thm:shift} gives the sharp value
$\Mfr(s)=2$ in the Beckermann--Reichel estimate~\eqref{eq:br}.
Indeed, the focal value is $2$, since
\[
 c+s=\frac{2s(1-\delta)}{1-k\delta},
 \qquad
 \sqrt{d^2+(c+s)^2}=\frac{(1+k)(1-\delta)}{1-k\delta}.
\]
For the scalar value, one has $C+s=s(1-\delta b)\geq0$, and hence
\[
 R+\abs{C+s}
 =
 1+s+b(1-s\delta)
 =
 1+s+\frac{(1-s)(1+s)\delta}{1+s\delta}
 \leq
 2,
\]
because $(1+s)\delta\leq1+s\delta$.

This computation does not use~\eqref{eq:br}. Neither does the proof of
Corollary~\ref{cor:upper} for polynomial symbols, which uses only
Proposition~\ref{prop:faber-dilation}, Theorem~\ref{thm:upper-annulus}
and Steps~1 and~2 of Lemma~\ref{lem:reduction}. After normalization and
a rotation as in Remark~\ref{rem:shift}{\rm(i)}, we therefore obtain an
independent proof of~\eqref{eq:br} on $K_\delta$ for polynomial symbols,
and with it the a priori bound of
Proposition~\ref{prop:boundary-calculus}{\rm(i)}. Theorem~\ref{thm:shift}
shows, in addition, that~\eqref{eq:br} is sharp for every prescribed value of
$g(0)$.

When both models attain the maximum, the shell construction will
require one symbol for both. To compare their phases, put
\[
 \mathcal N(t,m)=t+k\delta\overline t-s(1+\delta)m+s\bigl(1+(2-k)\delta\bigr),
 \qquad
 t\in\C,
 \quad
 m\in\R,
\]
and write
\[
 m_{\rm sc}=R+\abs{C+t},
 \qquad
 m_{\rm fo}=r+\sqrt{d^2+\abs{c+t}^2}
\]
for the two values in~\eqref{eq:M}.

\begin{lemma}\label{lem:phase}
Let $t\in\C$, and consider the scalar and focal models in the proof of
Theorem~\ref{thm:shift}.

\noindent{\rm (i)} If $t$ is real, both models can be built with
$\eta=1$.

\noindent{\rm (ii)} If $\Imag t\neq0$, then $\mathcal N(t,m)\neq0$ for every real
$m$, and
\[
 P_\omega Q_\omega=\frac{\mathcal N(t,m_{\rm sc})}{\overline{\mathcal N(t,m_{\rm sc})}}
 \quad\text{in the scalar model},
 \qquad
 v_1v_2=\frac{\mathcal N(t,m_{\rm fo})}{\overline{\mathcal N(t,m_{\rm fo})}}
 \quad\text{in the focal model}.
\]
In particular, if $m_{\rm sc}=m_{\rm fo}$, then the two models can be built
with the same $\eta$.
\end{lemma}

\begin{proof}
Two identities follow directly from~\eqref{eq:parameters}:
\begin{equation}\label{eq:phase-constants}
 \begin{aligned}
 (1+k\delta)C+s(1+\delta)R
 &=s(1+\delta)+sb(1-k\delta^2)
 =s\bigl(1+(2-k)\delta\bigr),\\
 (1+k\delta)c+s(1+\delta)r-s(1-\delta)d
 &=s(1+k\delta)+2s\delta a(1-k\delta)
 =s\bigl(1+(2-k)\delta\bigr).
 \end{aligned}
\end{equation}

For part~{\rm (i)}, let $t$ be real. Then $C+t$ and $c+t$ are real, and
in both models we may take $\omega=\pm1$. In the scalar model,
$P_\omega=Q_\omega=\omega$, so that $P_\omega Q_\omega=1$. In the focal
model, $v_1=v_2=\omega$ solve~\eqref{eq:phase}, so that $v_1v_2=1$. Thus
$\eta=1$ serves for both models.

For part~{\rm (ii)}, note first that
\[
 \Imag \mathcal N(t,m)=(1-k\delta)\Imag t\neq0.
\]

In the scalar model, $C+t\neq0$. Dividing the numerator and the
denominator of $P_\omega Q_\omega$ by $\omega$, and using
$\abs{\omega}=1$, we obtain
\[
 P_\omega Q_\omega=\frac{A}{\overline A},
 \qquad
 A=\omega+k\delta\overline\omega-s(1+\delta).
\]
Since $C$ is real, $\abs{C+t}\omega=C+t$ and
$\abs{C+t}\overline\omega=C+\overline t$. Using
$\abs{C+t}=m_{\rm sc}-R$ and the first identity
in~\eqref{eq:phase-constants},
\begin{align*}
 \abs{C+t}A
 &=
 t+k\delta\overline t+(1+k\delta)C+s(1+\delta)R-s(1+\delta)m_{\rm sc}\\
 &=
 \mathcal N(t,m_{\rm sc}).
\end{align*}
A positive factor does not change the quotient $A/\overline A$.

In the focal model, $\ell=\abs{c+t}>0$, so that $\mu,\nu>0$, and
$\omega=(c+t)/\ell$ is not real. Hence the coefficients
in~\eqref{eq:phase} do not vanish, and $v_1$, $v_2$ are determined.
Dividing the numerator and the denominator of $v_1v_2$ by $\omega\mu\nu$,
we obtain
\[
 v_1v_2=\frac{A'}{\overline{A'}},
 \qquad
 A'=\omega+k\delta\overline\omega-s\Bigl(\frac{\nu}{\mu}+\frac{\delta\mu}{\nu}\Bigr).
\]
Now $\ell\omega=c+t$, $\ell\overline\omega=c+\overline t$,
$\ell\nu/\mu=\kappa+d$, $\ell\mu/\nu=\kappa-d$ and $\kappa=m_{\rm fo}-r$.
By the second identity in~\eqref{eq:phase-constants},
\begin{align*}
 \ell A'
 &=
 t+k\delta\overline t+(1+k\delta)c+s(1+\delta)r-s(1-\delta)d
 -s(1+\delta)m_{\rm fo}\\
 &=
 \mathcal N(t,m_{\rm fo}).
\end{align*}
This proves part~{\rm (ii)}.

Finally, if $m_{\rm sc}=m_{\rm fo}$, parts~{\rm (i)} and~{\rm (ii)} show
that the two phase products coincide. Thus a common square root $\eta$
can be chosen.
\end{proof}

\begin{remark}\label{rem:shift}
\noindent{\rm (i)} If $g(0)=\lambda s$ with $\lambda\in\T$, linearity gives
\[
 \norm{f_g(T)+tI}
 =
 \norm{f_{\overline\lambda g}(T)+\overline\lambda tI}.
\]
Thus the supremum over pairs with $g(0)=\lambda s$ is
$\Mfr(\overline\lambda t)$, and the corresponding planar regions are
rotated by $\lambda$.

\noindent{\rm (ii)} The supremum in~\eqref{eq:main} does not change if $T$
is restricted to operators with $\overline{W(T)}\subset G_\delta$. Indeed,
the extremal symbols $g_\eta$ are holomorphic near $\overline\D$, so the
extremal operators may be replaced by $\theta T$ with $\theta\uparrow1$.

\noindent{\rm (iii)} For the focal model, one has $2\mu\nu=\ell/\kappa$, and
hence~\eqref{eq:focal-identity} gives
\begin{equation}\label{eq:focal-contact}
 \inner{f_{g_\eta}(T_v)x}{x}
 =
 c+2r\mu\nu\,\omega
 =
 c+\frac{r(c+t)}{\kappa}.
\end{equation}
This is the first coordinate of the focal contact point on the upper
Davis--Wielandt envelope.
\end{remark}

\subsection{The sharp bound on the quantum annulus}\label{ssec:annulus}

The upper bound of Subsection~\ref{ssec:upper} was proved on the quantum
annulus, and both extremal models live there. For matrices, the next
statement is equivalent to Theorem~\ref{thm:shift}. Indeed, Agler et
al.~\cite[Proposition~10.1]{AglerLykovaYoung} show that
$W(X+\delta X^{-1})\subseteq K_\delta$ for every such $X$, and
Theorem~\ref{thm:aly} gives the converse representation. The dilation
argument gives the estimate on every Hilbert space. Here the function
on the annulus is symmetric under $w\mapsto\delta/w$ and the value
$g(0)$ is prescribed. For the associated norm on the annulus, see~\cite{AglerLykovaYoungDp};
for spectral constants on the full annular function class,
see~\cite{Tsikalas,Pascoe,PalPascoeTomar}.

\begin{corollary}\label{cor:annulus}
Let $0<\delta<1$, $0\leq s<1$ and $t\in\C$. Let $X$ be an invertible
operator with
\[
 \norm{X}\leq1,
 \qquad
 \norm{\delta X^{-1}}\leq1,
\]
and let $g\in\Sch$ be holomorphic on a neighborhood of $\overline\D$. Then
\[
 \bigl\|g(X)+g(\delta X^{-1})+(t-s)I\bigr\|\leq M_{\delta,s}(t).
\]
Equality is attained, with $g=g_\eta$ for a suitable $\eta\in\T$, by a
unimodular scalar $X$ when the first value in~\eqref{eq:M} is the
larger one, and by a $2\times2$ matrix $X$ with $\norm{X}=1$ and
$X^2=\delta I$ otherwise.
\end{corollary}

\begin{proof}
Let $J$ and $U$ be given by Theorem~\ref{thm:mp}. Since $g$ is holomorphic
past $\T$, its Taylor series converges in norm at the four contractions
$X,\delta X^{-1},U,\delta U^{-1}$. Applying~\eqref{eq:laurent} term by term gives
\[
 g(X)+g(\delta X^{-1})+(t-s)I
 =
 J^*\bigl(g(U)+g(\delta U^{-1})+(t-s)I\bigr)J.
\]
The upper bound follows from Theorem~\ref{thm:upper-annulus}.
For the scalar value take $X=u$ from Lemma~\ref{lem:scalar-model}.
For the focal value take $X=U_0=A_v/\eta$ from
Lemma~\ref{lem:focal-model}. Then $\delta U_0^{-1}=U_0$ and
$2g_\eta(U_0)-sI=f_{g_\eta}(T_v)$. Theorem~\ref{thm:shift} gives
the asserted values.
\end{proof}

\section{The elliptic teardrop and sharp operator radii}\label{sec:consequences}

For $F=f_g(T)$, the shifted norms give support bounds on the first
moment $\inner{Fx}{x}$ and on the pair $(\inner{Fx}{x},\norm{Fx}^2)$.

\subsection{The elliptic teardrop}\label{ssec:teardrop}

We first determine the numerical-range region associated with the
shifted norm formula~\eqref{eq:main}. The scalar and focal
terms give the support functions of two disks. The relative position
of these disks is governed by
\[
 \tau=\frac{s(1+\delta)}{1+k\delta}.
\]

\begin{lemma}\label{lem:two-disks}
For $0<s<1$, one has
\[
 c-C
 =
 \frac{s(1-\delta)^2(1+k\delta)}
 {(1-k\delta)(1-k\delta^2)}
 >
 0,
 \qquad
 R-r=\tau(c-C),
 \qquad
 0<\tau<1.
\]
In particular, neither disk contains the other. For $s=0$, the two
disks coincide and have center $0$ and radius $1+\delta$.
\end{lemma}

\begin{proof}
Using~\eqref{eq:parameters} and clearing the common denominator gives
\[
 (c-C,R-r)
 =
 \frac{s(1-\delta)^2}
 {(1-k\delta)(1-k\delta^2)}
 \bigl(1+k\delta,s(1+\delta)\bigr).
\]
If $0<s<1$, then $c-C>0$, and
\[
 1+k\delta-s(1+\delta)
 =
 (1-s)(1-s\delta)
 >
 0.
\]
Thus $0<\tau<1$ and $\abs{R-r}<c-C$, which excludes containment in
either direction. At $s=0$, substitution in~\eqref{eq:parameters}
gives $C=c=0$ and $R=r=1+\delta$.
\end{proof}

\begin{theorem}\label{thm:teardrop}
For $0\leq\delta<1$ and $0\leq s<1$,
\begin{equation}\label{eq:hull}
 \bigcup_{(g,T)} W(f_g(T))
 =
 \mathcal T_{\delta,s}
 :=
 \conv\bigl(\overline D(C,R)\cup\overline D(c,r)\bigr),
\end{equation}
where the union is over all admissible pairs.
Every point of $\mathcal T_{\delta,s}$ is realized in dimension at
most three.

For $0<s<1$, the boundary consists of the arc of $\partial D(C,R)$
with outward normals $\zeta\in\T$ satisfying $\Real\zeta\leq\tau$, the
arc of $\partial D(c,r)$ with $\Real\zeta\geq\tau$, and the two common
tangent segments corresponding to $\Real\zeta=\tau$.
\end{theorem}

\begin{proof}
We first prove the inclusion from left to right. Let $(g,T)$ be
admissible, set $F=f_g(T)$, and let $z=\inner{Fx}{x}$ for a unit
vector $x$. For $\zeta\in\T$ and $\varrho>0$,
\begin{align*}
 2\varrho\Real(\overline\zeta z)
 &\leq
 \norm{(F+\varrho\zeta I)x}^2-\varrho^2\\
 &\leq
 M_{\delta,s}(\varrho\zeta)^2-\varrho^2.
\end{align*}
In~\eqref{eq:M}, expanding each square root at infinity gives
\[
 M_{\delta,s}(\varrho\zeta)
 =
 \varrho+\max\{R+C\Real\zeta,\ r+c\Real\zeta\}+O(\varrho^{-1}).
\]
Dividing the preceding inequality by $2\varrho$ and letting $\varrho\to\infty$ yields
\[
 \Real(\overline\zeta z)
 \leq
 \max\bigl\{R+C\Real\zeta,\ r+c\Real\zeta\bigr\}.
\]
The right-hand side is the support function of
$\mathcal T_{\delta,s}$. Thus $z$ belongs to that set.

For $s>0$, the difference of the two support values is
\[
 (r+c\Real\zeta)-(R+C\Real\zeta)
 =
 (c-C)(\Real\zeta-\tau).
\]
Thus the supporting line meets the hull at $C+R\zeta$ if
$\Real\zeta<\tau$, and at $c+r\zeta$ if $\Real\zeta>\tau$.
At $\Real\zeta=\tau$, a convex combination maximizes the support
functional exactly when each term with positive weight maximizes it.
The intersection with the supporting line is then the segment
joining $C+R\zeta$ and $c+r\zeta$. This proves the boundary
description. For $s=0$, the boundary is the common circle.

We next realize the boundary. Fix $\zeta\in\T$ and choose
$\eta\in\T$ with $\eta^2=P_\zeta Q_\zeta$.
Lemma~\ref{lem:scalar-model} gives $\lambda\in\partial K_\delta$ with
\[
 f_{g_\eta}(\lambda)=C+R\zeta.
\]
Applying
Lemma~\ref{lem:focal-model} with $\mu=\nu=1/\sqrt{2}$ gives a matrix
$T_\zeta$ unitarily equivalent to $A_\delta$ and a unit vector
$x_\zeta$, for the same symbol $g_\eta$, such that
\[
 \inner{f_{g_\eta}(T_\zeta)x_\zeta}{x_\zeta}
 =
 c+r\zeta.
\]
For $0\leq u\leq1$, put
\[
 T=\lambda\oplus T_\zeta,
 \qquad
 x=\sqrt{1-u}\oplus\sqrt{u}\,x_\zeta.
\]
Since $W(T)=K_\delta$, this is an admissible realization, and
\[
 \inner{f_{g_\eta}(T)x}{x}
 =
 (1-u)(C+R\zeta)+u(c+r\zeta).
\]
The endpoints and the segments at the two common supporting lines
therefore occur in dimension at most three.

Finally, the constant symbol $g=s$ realizes $s$, so the inclusion
already proved gives $s\in\mathcal T_{\delta,s}$. For
$z\in\mathcal T_{\delta,s}$ with $z\neq s$, the ray from $s$ through
$z$ meets the boundary at some $z_0$ beyond or equal to $z$, by
compactness and convexity. Write $z=(1-u)s+uz_0$, where $0<u\leq1$.
Choose a boundary realization
$z_0=\inner{f_{g_\eta}(T)x}{x}$ and set
\[
 g=(1-u)s+ug_\eta.
\]
Then $g\in\mathcal S_s\cap A(\overline\D)$, and linearity of the Faber
transform gives
\[
 \inner{f_g(T)x}{x}
 =
 (1-u)s+uz_0
 =
 z.
\]
This proves the reverse inclusion and the dimension bound.
\end{proof}

\begin{figure}[htbp]
\centering

\begin{minipage}{.485\textwidth}
\centering
\pgfmathsetmacro{\symb}{0.65}
\pgfmathsetmacro{\ecc}{0.1}
\pgfmathsetmacro{\kk}{\symb*\symb}
\pgfmathsetmacro{\scen}{-\symb*(1-\kk)*\ecc^2/(1-\kk*\ecc^2)}
\pgfmathsetmacro{\srad}{1+(1-\kk)*\ecc/(1-\kk*\ecc^2)}
\pgfmathsetmacro{\fcen}{\symb-2*\symb*(1-\kk)*\ecc/(1-\kk*\ecc)}
\pgfmathsetmacro{\frad}{(1-\kk)*(1+\ecc)/(1-\kk*\ecc)}
\pgfmathsetmacro{\ang}{acos(\symb*(1+\ecc)/(1+\kk*\ecc))}
\begin{tikzpicture}
\begin{axis}[
 width=\textwidth,height=.82\textwidth,
 axis equal image,
 xmin=-1.55,xmax=1.55,
 ymin=-1.42,ymax=1.42,
 xlabel={$\Real z$},
 ylabel={$\operatorname{Im}z$},
 title={(a) $\delta=0.1$},
 title style={font=\small},
 label style={font=\small},
 tick label style={font=\scriptsize},
 grid=major,
 grid style={gray!15}
]
\addplot[black,dashed,domain=0:360,samples=151]
 ({\scen+\srad*cos(x)},{\srad*sin(x)});

\addplot[red!70!black,densely dotted,domain=0:360,samples=151]
 ({\fcen+\frad*cos(x)},{\frad*sin(x)});

\addplot[black,thick,domain=\ang:360-\ang,samples=151]
 ({\scen+\srad*cos(x)},{\srad*sin(x)});

\addplot[red!70!black,thick,domain=-\ang:\ang,samples=81]
 ({\fcen+\frad*cos(x)},{\frad*sin(x)});

\addplot[blue!65!black,thick,domain=0:1,samples=2]
 ({(1-x)*(\scen+\srad*cos(\ang))+x*(\fcen+\frad*cos(\ang))},
  {((1-x)*\srad+x*\frad)*sin(\ang)});

\addplot[blue!65!black,thick,domain=0:1,samples=2]
 ({(1-x)*(\scen+\srad*cos(\ang))+x*(\fcen+\frad*cos(\ang))},
  {-((1-x)*\srad+x*\frad)*sin(\ang)});

\addplot[only marks,mark=*,mark size=1.8pt]
 coordinates {({\fcen+\frad},0)};
\end{axis}
\end{tikzpicture}
\end{minipage}\hfill
\begin{minipage}{.485\textwidth}
\centering
\pgfmathsetmacro{\symb}{0.65}
\pgfmathsetmacro{\ecc}{0.5}
\pgfmathsetmacro{\kk}{\symb*\symb}
\pgfmathsetmacro{\scen}{-\symb*(1-\kk)*\ecc^2/(1-\kk*\ecc^2)}
\pgfmathsetmacro{\srad}{1+(1-\kk)*\ecc/(1-\kk*\ecc^2)}
\pgfmathsetmacro{\fcen}{\symb-2*\symb*(1-\kk)*\ecc/(1-\kk*\ecc)}
\pgfmathsetmacro{\frad}{(1-\kk)*(1+\ecc)/(1-\kk*\ecc)}
\pgfmathsetmacro{\ang}{acos(\symb*(1+\ecc)/(1+\kk*\ecc))}
\begin{tikzpicture}
\begin{axis}[
 width=\textwidth,height=.82\textwidth,
 axis equal image,
 xmin=-1.55,xmax=1.55,
 ymin=-1.42,ymax=1.42,
 xlabel={$\Real z$},
 ylabel={$\operatorname{Im}z$},
 title={(b) $\delta=0.5$},
 title style={font=\small},
 label style={font=\small},
 tick label style={font=\scriptsize},
 grid=major,
 grid style={gray!15}
]
\addplot[black,dashed,domain=0:360,samples=151]
 ({\scen+\srad*cos(x)},{\srad*sin(x)});

\addplot[red!70!black,densely dotted,domain=0:360,samples=151]
 ({\fcen+\frad*cos(x)},{\frad*sin(x)});

\addplot[black,thick,domain=\ang:360-\ang,samples=151]
 ({\scen+\srad*cos(x)},{\srad*sin(x)});

\addplot[red!70!black,thick,domain=-\ang:\ang,samples=81]
 ({\fcen+\frad*cos(x)},{\frad*sin(x)});

\addplot[blue!65!black,thick,domain=0:1,samples=2]
 ({(1-x)*(\scen+\srad*cos(\ang))+x*(\fcen+\frad*cos(\ang))},
  {((1-x)*\srad+x*\frad)*sin(\ang)});

\addplot[blue!65!black,thick,domain=0:1,samples=2]
 ({(1-x)*(\scen+\srad*cos(\ang))+x*(\fcen+\frad*cos(\ang))},
  {-((1-x)*\srad+x*\frad)*sin(\ang)});

\addplot[only marks,mark=*,mark size=1.8pt]
 coordinates {({\scen-\srad},0)};
\end{axis}
\end{tikzpicture}
\end{minipage}

\caption{The elliptic teardrop $\mathcal T_{\delta,s}$ for $s=0.65$.
In both panels, the solid black and red arcs and the two blue segments
form its boundary. The dashed black circle is the scalar circle, and
the dotted red circle is the focal circle. In {\rm(a)}, the marked
point of largest modulus lies on the focal circle. In {\rm(b)}, it
lies on the scalar circle.}
\label{fig:arcs}
\end{figure}
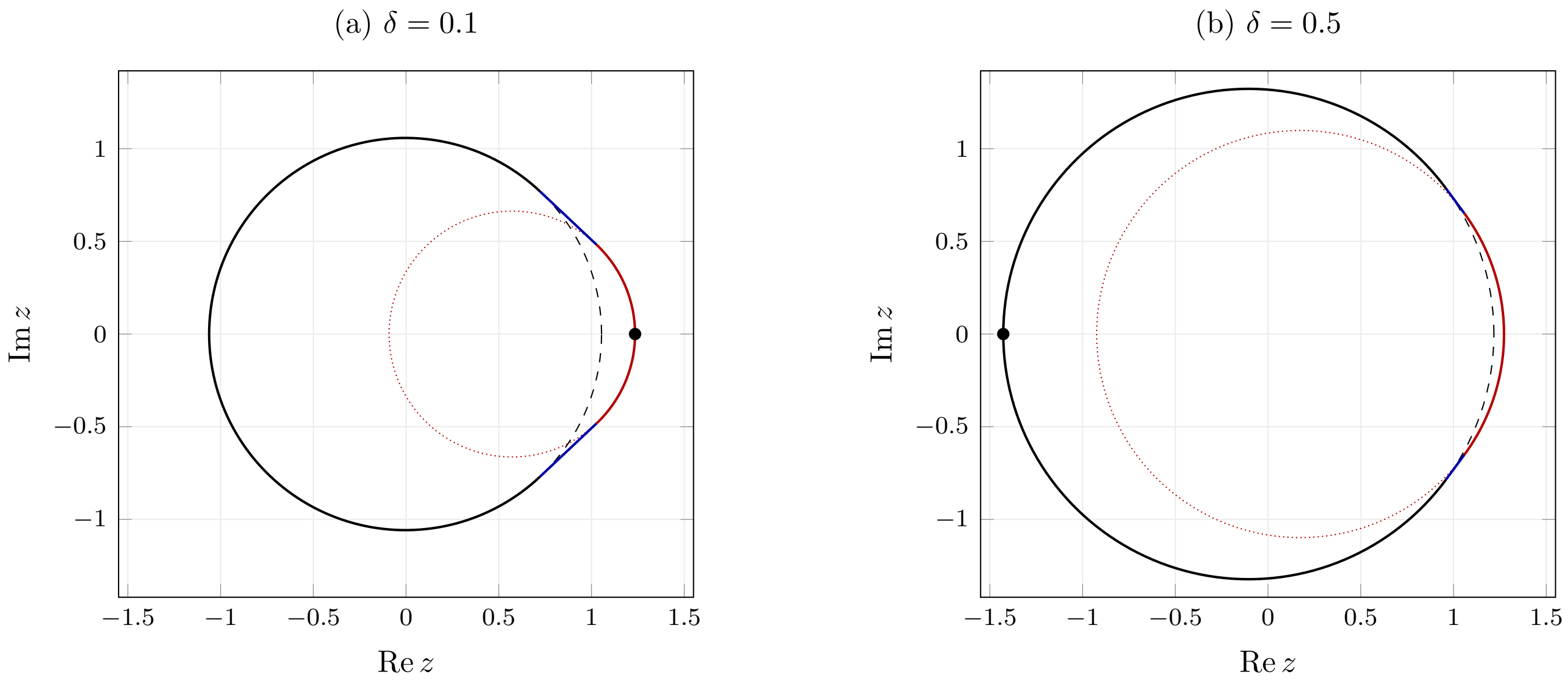

\begin{remark}\label{rem:teardrop}
\noindent{\rm (i)} For $s=0$, the teardrop is the disk of radius
$1+\delta$ centered at $0$, the smallest disk containing $K_\delta$. For
$\delta>0$, it contains the point $i(1+\delta)$, which lies outside
$K_\delta$. Thus, even for $g(0)=0$, the numerical range of $f_g(T)$ can
leave $K_\delta$, unlike the case $\delta=0$ of the Berger--Stampfli
theorem.

\noindent{\rm (ii)} Since $C\leq0$, the largest modulus on
$\mathcal T_{\delta,s}$ is
\[
 \max\{R-C,\ r+\abs{c}\}.
\]
The scalar disk gives the point $C-R$. The focal disk gives $c+r$ if
$c\geq0$ and $c-r$ if $c<0$. Figure~\ref{fig:arcs} shows that either
disk can contain the point of largest modulus.
\end{remark}

\FloatBarrier

\subsection{The upper Davis--Wielandt envelope}\label{ssec:envelope}

The shifted norm also controls the second coordinate of each
Davis--Wielandt pair.

Let $0\leq s<1$, and let $\mathcal A_{\delta,s}$ denote the union of
$\DW(f_g(T))$ over all admissible pairs. Put $M=M_{\delta,s}$. For each
$t\in\C$, define
\[
 \varphi_t(z,b)=b+2\Real(\overline t z),
 \qquad
 \Lambda(t)=M(t)^2-\abs{t}^2.
\]
For $F=f_g(T)$ and a unit vector $x$,
\begin{equation}\label{eq:phi-shell}
 \varphi_t\bigl(\inner{Fx}{x},\norm{Fx}^2\bigr)
 =
 \norm{(F+tI)x}^2-\abs{t}^2.
\end{equation}
Consequently, Theorem~\ref{thm:shift} gives
\[
 \sup_{(z,b)\in\mathcal A_{\delta,s}}\varphi_t(z,b)
 =
 \Lambda(t).
\]

The scalar and focal terms in~\eqref{eq:M} give rise to a disk
and a solid ellipsoid in $\C\times\mathbb R$. Define
\[
 \mathcal P
 =
 \left\{
 (C+v,C^2+R^2+2C\Real v):
 v\in\C,\ \abs{v}\leq R
 \right\}
\]
and
\[
 \mathcal E
 =
 \left\{
 (c+v,c^2+r^2+d^2+2c\Real v+2d\xi):
 (v,\xi)\in\C\times\mathbb R,\ \abs{v}^2+\xi^2\leq r^2
 \right\}.
\]
Set
\[
 \mathfrak K=\conv(\mathcal P\cup\mathcal E).
\]

\begin{lemma}\label{lem:body}
\noindent{\rm (i)} The first-coordinate projection of $\mathfrak K$ is
$\mathcal T_{\delta,s}$, and $b\geq\abs{z}^2$ for every
$(z,b)\in\mathfrak K$.

\noindent{\rm (ii)} For $t\in\C$, put
$\kappa(t)=\sqrt{d^2+\abs{c+t}^2}$. Then
\begin{align*}
 \max_{(z,b)\in\mathcal P}\varphi_t(z,b)
 &=
 (R+\abs{C+t})^2-\abs{t}^2,\\
 \max_{(z,b)\in\mathcal E}\varphi_t(z,b)
 &=
 (r+\kappa(t))^2-\abs{t}^2.
\end{align*}
If $C+t\neq0$, the maximum on $\mathcal P$ is attained only at
$v=R(C+t)/\abs{C+t}$ in the parametrization of $\mathcal P$.
The maximum on $\mathcal E$ is attained only at
$(v,\xi)=r(c+t,d)/\kappa(t)$.

\noindent{\rm (iii)} The maximum on $\mathfrak K$ equals $\Lambda(t)$.
Its maximizing set is the convex hull of the maximizing points from
those of $\mathcal P$ and $\mathcal E$ whose maximum equals
$\Lambda(t)$.
\end{lemma}

\begin{proof}
For part~{\rm(i)}, the projections of $\mathcal P$ and $\mathcal E$
are the two disks in~\eqref{eq:hull}. On $\mathcal P$,
\[
 b-\abs{z}^2=R^2-\abs{v}^2\geq0.
\]
On $\mathcal E$, we have
\begin{align*}
 b-\abs{z}^2
 &=
 r^2-\abs{v}^2+d^2+2d\xi\\
 &\geq
 \bigl(\sqrt{r^2-\abs{v}^2}-d\bigr)^2
 \geq0,
\end{align*}
because $\xi\geq-\sqrt{r^2-\abs{v}^2}$ and $d>0$. Convexity of the
epigraph of $z\mapsto\abs{z}^2$ extends the inequality to $\mathfrak K$.

For part~{\rm(ii)}, use the coordinates in the definitions of
$\mathcal P$ and $\mathcal E$. The quantity
$\varphi_t(z,b)+\abs{t}^2$ takes the forms
\[
 \begin{cases}
 \abs{C+t}^2+R^2+2\Real\bigl(\overline{(C+t)}v\bigr),
 & (z,b)\in\mathcal P,\\[3pt]
 \kappa(t)^2+r^2+2\Real\bigl(\overline{(c+t)}v\bigr)+2d\xi,
 & (z,b)\in\mathcal E.
 \end{cases}
\]
Cauchy--Schwarz gives the two maxima and the stated equality
conditions. The vector $(c+t,d)$ is nonzero because $d>0$.

The maximum over the convex hull is the larger of the two maxima,
which is $\Lambda(t)$ by~\eqref{eq:M}. A convex combination
maximizes $\varphi_t$ exactly when every term with positive weight
maximizes it. This proves part~{\rm(iii)}.
\end{proof}

For $z\in\mathcal T_{\delta,s}$, define
\[
 b^*(z)
 =
 \max\{b:(z,b)\in\mathfrak K\}.
\]
Since $\mathfrak K$ is compact and convex, $b^*$ is upper
semicontinuous and concave on $\mathcal T_{\delta,s}$.

\begin{theorem}\label{thm:dw}
Let $0\leq\delta<1$ and $0\leq s<1$.

\noindent{\rm (i)} For every $z\in\mathcal T_{\delta,s}$,
\begin{equation}\label{eq:envelope}
 b^*(z)
 =
 \inf_{t\in\C}\bigl\{\Lambda(t)-2\Real(\overline t z)\bigr\}
 =
 \max\{b:(z,b)\in\mathcal A_{\delta,s}\}.
\end{equation}
Every upper pair $(z,b^*(z))$ is realized by an automorphism $g_\eta$
and an operator unitarily equivalent to $A_\delta\oplus\lambda$ on
$\C^3$, with $\lambda\in\partial K_\delta$.

\noindent{\rm (ii)} For $1\leq\rho\leq2$, the smallest closed convex set
containing $W^\rho(f_g(T))$ for all admissible pairs is
\[
 \overline{\conv}\Phi_\rho(\mathfrak K).
\]
\end{theorem}

\begin{proof}
\noindent\emph{Step 1: the envelope formula.}
The set
\[
 \mathfrak K^\downarrow
 =
 \bigl\{(z,b):z\in\mathcal T_{\delta,s},\ b\leq b^*(z)\bigr\}
\]
is closed and convex by the properties of $b^*$. The maximum of
$\varphi_t$ over this set is the same as over $\mathfrak K$, namely
$\Lambda(t)$. Thus
\[
 b^*(z)
 \leq
 \Lambda(t)-2\Real(\overline t z),
 \qquad t\in\C.
\]

Conversely, take $z\in\mathcal T_{\delta,s}$ and $b>b^*(z)$.
A linear functional strictly separating $(z,b)$ from
$\mathfrak K^\downarrow$ has a nonnegative coefficient of the second
coordinate, since $\mathfrak K^\downarrow$ contains the downward
vertical ray from every point of $\mathfrak K$. This coefficient
must be positive, since $(z,b^*(z))\in\mathfrak K^\downarrow$ has the
same first coordinate as $(z,b)$. Normalizing the functional gives
$\varphi_t$ for some $t\in\C$, and hence
\[
 b+2\Real(\overline t z)>\Lambda(t).
\]
This gives the first equality in~\eqref{eq:envelope} by letting
$b\downarrow b^*(z)$.

If $(z,b)\in\mathcal A_{\delta,s}$, then
$z\in\mathcal T_{\delta,s}$ by Theorem~\ref{thm:teardrop}, and
$\varphi_t(z,b)\leq\Lambda(t)$ for every $t$ by~\eqref{eq:phi-shell}.
The infimum formula therefore gives $b\leq b^*(z)$.

\smallskip
\noindent\emph{Step 2: realization on the supporting planes.}
Fix $t\in\C$. If the scalar term in~\eqref{eq:M} attains $M(t)$,
then $C+t\neq0$. Indeed, at $t=-C$ the focal value satisfies
\[
 r+\sqrt{d^2+(c-C)^2}
 >
 r+c-C
 \geq
 R,
\]
by Lemma~\ref{lem:two-disks}. By Lemma~\ref{lem:body}{\rm(ii)}, the
scalar maximizing pair is
\[
 \bigl(z_{\rm sc},\abs{z_{\rm sc}}^2\bigr),
 \qquad
 z_{\rm sc}
 =
 C+R\frac{C+t}{\abs{C+t}}.
\]
Lemma~\ref{lem:scalar-model} realizes it at a boundary scalar
$\lambda$.

If the focal term attains $M(t)$, the maximizing pair on
$\mathcal E$ has first coordinate
\[
 z_{\rm fo}=c+\frac{r(c+t)}{\kappa(t)}.
\]
For the focal model of Theorem~\ref{thm:shift},
\eqref{eq:focal-contact} gives this first coordinate, and
\eqref{eq:phi-shell} gives $\Lambda(t)$ as the value of
$\varphi_t$. The first coordinate and the value of $\varphi_t$
determine the second coordinate, so the model realizes the
maximizing pair on $\mathcal E$.

If both terms attain $M(t)$, Lemma~\ref{lem:phase} gives a common
symbol $g_\eta$ for these models. Write $T_v$ for the focal matrix
in this realization and $x$ for its unit vector. For $0\leq u\leq1$,
the pair for $\lambda\oplus T_v$ at
\[
 \sqrt{1-u}\oplus\sqrt{u}\,x
\]
is the convex combination of the two maximizing pairs with weights
$1-u$ and $u$, since both moments add with these weights along the
direct sum. By Lemma~\ref{lem:body}{\rm(iii)}, this realizes the
whole maximizing set. If only one term attains $M(t)$, the
corresponding realization can be made on the same three-dimensional
model by adjoining $A_\delta$ to the scalar or any boundary scalar
to the focal matrix. The vector is extended by zero.

\smallskip
\noindent\emph{Step 3: realization of every upper pair.}
Suppose first that $z$ belongs to the interior of
$\mathcal T_{\delta,s}$. At $(z,b^*(z))$, a supporting functional of
$\mathfrak K^\downarrow$ has a nonnegative coefficient of $b$ because
of the downward rays. If that coefficient were zero, a nonzero
linear functional of $z$ would attain its maximum over
$\mathcal T_{\delta,s}$ at an interior point. Thus the coefficient
is positive, and after normalization $(z,b^*(z))$ maximizes some
$\varphi_t$. Step~2 realizes this pair.

Now let $z_0\in\partial\mathcal T_{\delta,s}$ and fix an interior
point $z_1$. For $0<u<1$, the point $z_u=(1-u)z_0+uz_1$ is interior.
Concavity gives
\[
 b^*(z_u)\geq(1-u)b^*(z_0)+ub^*(z_1).
\]
Together with upper semicontinuity, this implies
$b^*(z_u)\to b^*(z_0)$ as $u\downarrow0$.

The pairs for $f_{g_\eta}(A_\delta\oplus\lambda)$ at unit vectors in
$\C^3$ form a compact set. Indeed, the parameters $\eta$, $\lambda$
and the unit vector vary over compact sets, and the moments depend
continuously on them by the rational formula~\eqref{eq:mob-lift}.
This set contains $(z_u,b^*(z_u))$ for every $0<u<1$, after absorbing
the unitary equivalence of the focal model into the vector.
It therefore contains the limit $(z_0,b^*(z_0))$. This proves
part~{\rm(i)}.

\smallskip
\noindent\emph{Step 4: deformed numerical ranges.}
Let
\[
 \mathcal D
 =
 \overline{\conv}\Phi_\rho(\mathcal A_{\delta,s}).
\]
By definition, $\mathcal D$ is the smallest closed convex set
containing all the deformed numerical ranges in part~{\rm(ii)}.
For every $(z,b)$ in $\mathcal A_{\delta,s}$ or in $\mathfrak K$,
\[
 \abs{z}^2\leq b\leq b^*(z).
\]
As the second coordinate increases, $\Phi_\rho$ moves outward on
the ray determined by $z$ when $z\neq0$, and it equals zero when
$z=0$. Thus $\Phi_\rho(z,b)$ lies on the segment from $0$ to
$\Phi_\rho(z,b^*(z))$.

The circle $C+R\T$ at height $b=\abs{z}^2$ belongs to $\mathcal P$
and is realized by Lemma~\ref{lem:scalar-model}. At this height
$\Phi_\rho(z,b)=z$. Since
\[
 R-\abs{C}
 =
 1+(1-s\delta)(R-1)
 >
 0,
\]
both closed convex hulls contain $0$. By part~{\rm(i)}, both also
contain $\Phi_\rho(z,b^*(z))$ for every $z\in\mathcal T_{\delta,s}$.
Thus the two hulls coincide, proving part~{\rm(ii)}.
\end{proof}

\subsection{The sharp operator-radius bound}\label{ssec:radius}

For every $1\leq\rho\leq2$, a single symbol $\beta_s$ attains the sharp
output bound. The two possible values follow by inserting the shifted
norm formula into the positivity criterion for $w_\rho$.

\begin{corollary}\label{cor:radius}
For $0\leq\delta<1$, $0\leq s<1$ and $1\leq\rho\leq2$, one has
\begin{equation}\label{eq:radius}
 \Gamma_{\delta,\rho}(s)
 :=
 \sup_{g\in\mathcal S_s}\norm{f_g}_{\rho,\mathrm{bfd}}
 =
 \max\left\{
 R-C,\
 \frac{r+p\abs{c}
 +\sqrt{(r+p\abs{c})^2-\rho q(r^2-d^2-c^2)}}{\rho}
 \right\}.
\end{equation}
The same value is the supremum of $w_\rho(f_g(T))$ over all admissible
pairs. For the symbol $g_1=\beta_s$, the two values are realized on the
closed ellipse by $T=-(1+\delta)$ and $T=A_\delta$, respectively.

Consequently, if $g\in H^\infty(\D)$ is nonzero and
$m=\norm{g}_\infty$, then
\[
 \norm{f_g}_{\rho,\mathrm{bfd}}
 \leq
 m\,\Gamma_{\delta,\rho}\left(\frac{\abs{g(0)}}{m}\right),
 \qquad
 \Gamma_{\delta,\rho}(1)=1.
\]
\end{corollary}

\begin{proof}
Let $F=f_g(T)$ for an admissible pair $(g,T)$. Suppose first
that $1\leq\rho<2$, so $q>0$. For $h>0$ and
$\zeta\in\T$, the identity $p^2+\rho q=1$ gives
\begin{align*}
 h^2I-(qF+ph\zeta I)^*(qF+ph\zeta I)
 &=
 q\bigl(\rho h^2I-2ph\Real(\overline\zeta F)-qF^*F\bigr).
\end{align*}
By Lemma~\ref{lem:pencil}{\rm(i)}, $w_\rho(F)\leq h$ is therefore
equivalent to
\[
 \norm{qF+ph\zeta I}\leq h
 \qquad(\zeta\in\T).
\]
By Theorem~\ref{thm:shift} and Remark~\ref{rem:shift}{\rm(ii)}, the
condition holds for all $F=f_g(T)$ exactly when
\[
 qM\Bigl(\frac{ph}{q}\zeta\Bigr)\leq h,
 \qquad
 \zeta\in\T.
\]
Maximizing the two terms separately over $\zeta$ gives
\[
 q(R-C)+ph\leq h,
 \qquad
 qr+\sqrt{q^2d^2+(q\abs{c}+ph)^2}\leq h.
\]
The first condition is $h\geq R-C$. The second is equivalent to
$h\geq qr$ together with
\[
 \rho h^2-2(r+p\abs{c})h+q(r^2-d^2-c^2)\geq0.
\]
At $h=qr$, the polynomial equals
\[
 -q\bigl((pr+\abs{c})^2+d^2\bigr)<0.
\]
Thus, for $h\geq qr$, the inequality holds precisely from the larger
root onward. This root
is the second value in~\eqref{eq:radius}.

For $\rho=2$, \eqref{eq:ppw} and Theorem~\ref{thm:teardrop} reduce
the upper bound to the largest modulus on $\mathcal T_{\delta,s}$.
By Remark~\ref{rem:teardrop}{\rm(ii)}, that is
\[
 \max\{R-C,\ r+\abs{c}\},
\]
which agrees with~\eqref{eq:radius} at $p=1$, $q=0$.

We next prove sharpness for the single symbol $g_1$. By
Lemma~\ref{lem:scalar-model},
\[
 f_{g_1}\bigl(-(1+\delta)\bigr)=C-R.
\]
For the focal model, let $F_1=f_{g_1}(A_\delta)=cI+aA_\delta$.
The focal bound in the proof of Theorem~\ref{thm:shift} applies to
$g_1$ for every complex shift, since the spectrum of $\Sigma$ on this
model consists of the two endpoints $\pm j_0$ when $\delta>0$. For
$\delta=0$ the nilpotent computation in the same proof applies.
For real $t$, Lemma~\ref{lem:phase}{\rm(i)} gives $\eta=1$ in the
attaining model, and $T_v$ is unitarily equivalent to $A_\delta$. Thus
\[
 \norm{F_1+tI}\leq r+\sqrt{d^2+\abs{c+t}^2},
 \qquad\text{with equality for }t\in\R.
\]

When $q>0$, this estimate applied at $t=ph\zeta/q$ yields
\[
 \max_{\zeta\in\T}\norm{qF_1+ph\zeta I}
 =
 qr+\sqrt{q^2d^2+(q\abs{c}+ph)^2}.
\]
Equality follows by taking $\zeta=1$ when $c\geq0$ and $\zeta=-1$
when $c<0$. The same norm criterion therefore gives exactly the
focal value in~\eqref{eq:radius}. When $\rho=2$, the ellipse
$W(F_1)=c+aK_\delta$ has real center $c$ and semi-axes $r\geq d$,
so its largest modulus is $r+\abs{c}$.

The transform $f_{g_1}$ is holomorphic on a neighborhood of
$K_\delta$, so Lemma~\ref{lem:closed-class} transfers the values
from the two models to $\norm{f_{g_1}}_{\rho,\mathrm{bfd}}$. Finally, the normalized
bound follows by applying~\eqref{eq:radius} to $g/m$ after the
rotation from Remark~\ref{rem:shift}{\rm(i)}, and using
$f_{g/m}=f_g/m$. If $\abs{g(0)}=m$, the symbol is constant by the
maximum modulus principle, which gives the endpoint $s=1$.
\end{proof}

\begin{remark}\label{rem:radius}
\noindent{\rm (i)} At $\rho=1$, \eqref{eq:radius} gives the shifted norm
formula at $t=0$. At $\rho=2$, it gives the largest modulus on
$\mathcal T_{\delta,s}$.

\noindent{\rm (ii)} If $s=0$, then $C=c=0$, $R=r=1+\delta$ and
$d=1-\delta$. The focal model is $A_\delta$ itself. Since
$w_\rho(A_\delta)\geq w(A_\delta)=1+\delta$, the focal value dominates
the scalar one. Thus
\[
 \Gamma_{\delta,\rho}(0)
 =
 w_\rho(A_\delta)
 =
 \frac{1+\delta+\sqrt{(1+\delta)^2-4\rho(2-\rho)\delta}}{\rho}.
\]

\noindent{\rm (iii)} If $\delta=0$, the focal value is at least
$w(F_1)=1+s-s^2\geq1$, while the scalar value is $1$. Therefore
\[
 \Gamma_{0,\rho}(s)
 =
 \frac{1-s^2+ps+
 \sqrt{(1-s^2+ps)^2+\rho qs^2}}{\rho}.
\]
At $\rho=2$, this is Drury's bound $1+s-s^2$ for the numerical
radius~\cite{Drury}. At $\rho=1$, it is the sharp norm bound for
$w(T)\leq1$ of~\cite{Drury,SchwenningerDeVries}.
\end{remark}

\section{Equality and rigidity}\label{sec:equality}

\subsection{Equality by contact}\label{ssec:contact}

For a fixed shift $t$, identity~\eqref{eq:phi-shell} expresses norm
attainment as contact with a supporting plane of the upper
Davis--Wielandt envelope. The two terms in~\eqref{eq:M} determine
where contact occurs.

Throughout this subsection, let $0\leq\delta<1$ and $0\leq s<1$, and put
$M=M_{\delta,s}$. We regard the two values $m_{\rm sc}$ and
$m_{\rm fo}$ of Subsection~\ref{ssec:attainment} as functions of $t$. Thus
$m_{\rm fo}(t)=r+\kappa(t)$, with $\kappa$ as in
Lemma~\ref{lem:body}{\rm(ii)}, and $M(t)=\max\{m_{\rm sc}(t),m_{\rm fo}(t)\}$. If
$m_{\rm sc}(t)=M(t)$, then $C+t\neq0$ by Step~2 of the proof of
Theorem~\ref{thm:dw}, and we put
\[
 z_{\rm sc}(t)=C+R\,\frac{C+t}{\abs{C+t}}.
\]
Since $\kappa(t)\geq d>0$, we also put
\[
 z_{\rm fo}(t)=c+r\,\frac{c+t}{\kappa(t)}.
\]
Let $Z(t)\subset\C$ be the convex hull of $z_{\rm sc}(t)$, included when
$m_{\rm sc}(t)=M(t)$, and of $z_{\rm fo}(t)$, included when
$m_{\rm fo}(t)=M(t)$. With $\Lambda(t)=M(t)^2-\abs{t}^2$ as in
Subsection~\ref{ssec:envelope}, define
\[
 \mathcal G(t)
 =
 \bigl\{\bigl(z,\Lambda(t)-2\Real(\overline tz)\bigr):z\in Z(t)\bigr\}.
\]
Lemma~\ref{lem:body} identifies $\mathcal G(t)$ with the exposed face
of $\mathfrak K$ on the plane $\varphi_t=\Lambda(t)$.

\begin{theorem}\label{thm:equality-contact}
Let $0\leq\delta<1$ and $0\leq s<1$. Let $(g,T)$ be admissible, and put
$F=f_g(T)$.

\noindent{\rm (i)} For every $t\in\C$,
\[
 \norm{F+tI}=M(t)
 \quad\Longleftrightarrow\quad
 \overline{\DW(F)}\cap\mathcal G(t)\neq\varnothing.
\]
A unit vector $x$ satisfies $\norm{(F+tI)x}=M(t)$ if and only if
\[
 \bigl(\inner{Fx}{x},\norm{Fx}^2\bigr)\in\mathcal G(t).
\]

\noindent{\rm (ii)} Let $1\leq\rho<2$, put $h=\Gamma_{\delta,\rho}(s)$,
and write
\[
 t_\zeta=\frac{ph}{q}\zeta,
 \qquad
 \zeta\in\T.
\]
Then $w_\rho(F)=h$ if and only if there is $\zeta\in\T$ such that
\[
 qM(t_\zeta)=h
 \qquad\text{and}\qquad
 \overline{\DW(F)}\cap\mathcal G(t_\zeta)\neq\varnothing.
\]
For $\rho=2$, one has $w(F)=\Gamma_{\delta,2}(s)$ if and only if
$\overline{W(F)}$ contains a point of modulus $\max\{R-C,\ r+\abs{c}\}$.

\noindent If $\mathcal H$ is finite-dimensional, the closures may be
omitted.
\end{theorem}

\begin{proof}
By~\eqref{eq:phi-shell},
\[
 \norm{F+tI}^2-\abs{t}^2
 =
 \max_{(z,b)\in\overline{\DW(F)}}\varphi_t(z,b).
\]
The maximum exists because $\overline{\DW(F)}$ is compact.
By Theorem~\ref{thm:dw}{\rm(i)} and upper semicontinuity of $b^*$,
every point of this closure satisfies $b\leq b^*(z)$. Since $\varphi_t$
increases strictly with $b$, it can equal $\Lambda(t)$ only at an upper
pair on $\mathcal G(t)$. Conversely every pair on that face gives
equality. Applying the same identity at a unit vector proves the
assertion about attainment.

For part~{\rm (ii)}, let $q>0$. By the identity in the proof of
Corollary~\ref{cor:radius} and Lemma~\ref{lem:pencil}{\rm(i)}, for every
$h_0>0$,
\[
 w_\rho(F)\leq h_0
 \quad\Longleftrightarrow\quad
 \psi(h_0):=\max_{\zeta\in\T}\norm{qF+ph_0\zeta I}\leq h_0.
\]
Theorem~\ref{thm:shift} and the proof of Corollary~\ref{cor:radius} give
$w_\rho(F)\leq h$, and hence $\psi(h)\leq h$.

We claim that $w_\rho(F)=h$ if and only if $\psi(h)=h$. If $\psi(h)<h$,
the continuity of $\psi$ gives $\psi(h_0)<h_0$ for some $h_0<h$, and hence
$w_\rho(F)<h$. Conversely, if $w_\rho(F)\leq h_0<h$, then
\[
 \psi(h)
 \leq
 \psi(h_0)+p(h-h_0)
 \leq
 h_0+p(h-h_0)
 <
 h,
\]
since $p<1$.

Now let $\psi(h)=h$. By compactness of $\T$, there is $\zeta\in\T$ with
$\norm{qF+ph\zeta I}=h$. By Theorem~\ref{thm:shift} and the proof of
Corollary~\ref{cor:radius},
\[
 h
 =
 q\norm{F+t_\zeta I}
 \leq
 qM(t_\zeta)
 \leq
 h.
\]
Hence $qM(t_\zeta)=h$ and $\norm{F+t_\zeta I}=M(t_\zeta)$, and part~{\rm (i)}
applies. Conversely, under the two stated conditions, part~{\rm (i)}
gives
\[
 \norm{qF+ph\zeta I}=qM(t_\zeta)=h,
\]
so that $\psi(h)=h$.

For $\rho=2$, one has $w_2=w$, and Corollary~\ref{cor:radius} gives
$\Gamma_{\delta,2}(s)=\max\{R-C,\ r+\abs{c}\}$. Since $\overline{W(F)}$
is compact, $w(F)$ equals this value if and only if $\overline{W(F)}$
contains a point of this modulus.

Finally, if $\mathcal H$ is finite-dimensional, then $\DW(F)$ and $W(F)$
are compact.
\end{proof}

On a strict scalar branch, contact has a further consequence.

\begin{remark}\label{rem:scalar-contact}
Suppose that $m_{\rm sc}(t)>m_{\rm fo}(t)$, and let $x$ be a unit vector
with $\norm{(F+tI)x}=M(t)$. Then $\C x$ reduces $F$, and
\[
 Fx=z_{\rm sc}(t)x.
\]
Indeed, in this case $\mathcal G(t)$ consists of the single point
$(z_{\rm sc}(t),\abs{z_{\rm sc}(t)}^2)$, and
Theorem~\ref{thm:equality-contact}{\rm(i)} places the pair of $x$ at this
point. Hence
$\norm{Fx}^2=\abs{\inner{Fx}{x}}^2$, and the Cauchy--Schwarz inequality
gives $Fx=zx$ with $z=z_{\rm sc}(t)$. Since $x$ attains the norm of
$F+tI$,
\[
 (F+tI)^*(F+tI)x=M(t)^2x.
\]
As $(F+tI)x=(z+t)x$ and $\abs{z+t}=M(t)>0$, it follows that
$F^*x=\overline zx$.
\end{remark}

\subsection{Davis--Wielandt shells of the elliptical class}\label{ssec:shells}

The extremizing automorphisms have linear-fractional Faber transforms
by~\eqref{eq:mob-lift}. We now fix such a function $f$ and ask which
operators attain its sharp bound. The contact in
Theorem~\ref{thm:equality-contact} concerns the output $f(T)$. For a
linear-fractional function the denominator allows us to express the
same pencil in terms of the two moments of the input $T$.

The proof uses the dilation theorem for $2\times2$ numerical ranges, due
to Choi and Li~\cite{ChoiLi}; see~\cite[Theorem~1.2]{LiPoon} for a unified
proof. For $\delta=0$, the upper estimate in Proposition~\ref{prop:shell}
is the local estimate of Klaja et al.~\cite{KlajaMashreghiRansford}.

\begin{theorem}[\cite{ChoiLi,LiPoon}]\label{thm:choi-li}
If $W(T)\subseteq K_\delta=W(A_\delta)$, then there are a Hilbert space
$\mathcal M$ and an isometry $J:\mathcal H\to\C^2\otimes\mathcal M$ such
that
\[
 T=J^*(A_\delta\otimes I_{\mathcal M})J.
\]
\end{theorem}

\begin{proposition}\label{prop:shell}
Let $0\leq\delta<1$, and put
\[
 \mathfrak b_\delta(z)
 =
 2(1+\delta^2)+2(1-\delta^2)\sqrt{1-m_\delta(z)^2},
 \qquad
 z\in K_\delta.
\]
Then
\[
 \bigcup_{W(T)\subseteq K_\delta}\DW(T)
 =
 \bigl\{(z,b):z\in K_\delta,\ \abs{z}^2\leq b\leq\mathfrak b_\delta(z)\bigr\}.
\]
On each fiber, the lower endpoint is attained by the scalar $z$ and the
upper endpoint by $A_\delta$.
\end{proposition}

\begin{proof}
Let $W(T)\subseteq K_\delta$, and let $x$ be a unit vector. By
Theorem~\ref{thm:choi-li},
\[
 T=J^*(A_\delta\otimes I_{\mathcal M})J
\]
for an isometry $J$. Write $Jx=(\xi_1,\xi_2)\in\mathcal M\oplus\mathcal M$,
and put
\[
 \chi=2\inner{\xi_2}{\xi_1},
 \qquad
 \vartheta=\norm{\xi_2}^2,
\]
so that $\norm{\xi_1}^2=1-\vartheta$. As in Lemma~\ref{lem:focal},
\[
 z:=\inner{Tx}{x}=\chi+\delta\overline\chi,
 \qquad
 \abs{\chi}=m_\delta(z),
 \qquad
 \abs{\chi}^2\leq4\vartheta(1-\vartheta).
\]
Since $J^*$ is a contraction,
\[
 \norm{Tx}^2
 \leq
 \norm{(A_\delta\otimes I)Jx}^2
 =
 4\vartheta+4\delta^2(1-\vartheta).
\]

The inequality $\abs{\chi}^2\leq4\vartheta(1-\vartheta)$ gives
\[
 \vartheta\leq \vartheta_+:=\frac{1+\sqrt{1-\abs{\chi}^2}}{2},
\]
and the bound for $\norm{Tx}^2$ increases with $\vartheta$. Substituting $\vartheta_+$, we
obtain
\[
 \norm{Tx}^2\leq\mathfrak b_\delta(z).
\]
Moreover, the Cauchy--Schwarz inequality gives $\abs{z}^2\leq\norm{Tx}^2$.

Conversely, fix $z\in K_\delta$, and put $\chi=L_\delta^{-1}(z)$, so that
$\abs{\chi}=m_\delta(z)\leq1$. With $\vartheta=\vartheta_+\geq1/2$, the vector
\[
 y=\Bigl(\frac{\overline\chi}{2\sqrt\vartheta},\sqrt\vartheta\Bigr)
\]
is a unit vector, because $4\vartheta(1-\vartheta)=\abs{\chi}^2$. Lemma~\ref{lem:focal}
gives
\[
 \inner{A_\delta y}{y}=z,
 \qquad
 \norm{A_\delta y}^2=4\vartheta+4\delta^2(1-\vartheta)=\mathfrak b_\delta(z).
\]
The scalar $z$ gives the pair $(z,\abs{z}^2)$. Moreover,
$W(A_\delta\oplus z)=K_\delta$. For $0\leq u\leq1$, the unit vector
$\sqrt u\,y\oplus\sqrt{1-u}$ gives the pair
\[
 \bigl(z,\ u\mathfrak b_\delta(z)+(1-u)\abs{z}^2\bigr)
\]
for $A_\delta\oplus z$.
\end{proof}

\subsection{The norm and the numerical range of \texorpdfstring{$f(T)$}{f(T)}}\label{ssec:lf}

Proposition~\ref{prop:shell} reduces every inequality that is affine in
the two moments $\inner{Tx}{x}$ and $\norm{Tx}^2$ to scalars and to the
focal matrix. For a linear-fractional function $f$, both the
positivity criterion for $w_\rho(f(T))$ and the support functions of
$W(f(T))$ lead to such inequalities after congruence by the denominator.

\begin{theorem}\label{thm:lf}
Let
\[
 f(z)=\frac{\alpha z+\alpha'}{\gamma z+\gamma'},
 \qquad
 \gamma z+\gamma'\neq0\quad\text{on }K_\delta,
\]
and let $1\leq\rho\leq2$. Then
\begin{align}
 \norm{f}_{\rho,\mathrm{bfd}}
 &=
 \max\bigl\{\norm{f}_{\infty,K_\delta},\ w_\rho(f(A_\delta))\bigr\},
 \label{eq:individual}\\
 \bigcup_{W(T)\subseteq K_\delta}W(f(T))
 &=
 \conv\bigl(f(K_\delta)\cup W(f(A_\delta))\bigr).
 \label{eq:individual-range}
\end{align}
Moreover,
\[
 f(A_\delta)=a_fI+b_fA_\delta,
 \qquad
 W(f(A_\delta))=a_f+b_fK_\delta,
\]
where
\[
 a_f=\frac{f(2\sqrt\delta)+f(-2\sqrt\delta)}{2},
 \qquad
 b_f=\frac{f(2\sqrt\delta)-f(-2\sqrt\delta)}{4\sqrt\delta}
 \qquad
 (\delta>0),
\]
and $a_f=f(0)$, $b_f=f'(0)$ for $\delta=0$.
\end{theorem}

\begin{proof}
\emph{Step 1: the focal data.} Let $\delta>0$, and put
\[
 \varsigma=2\sqrt\delta,
 \qquad
 P_\pm=\frac12\Bigl(I\pm\frac{A_\delta}{\varsigma}\Bigr).
\]
By Lemma~\ref{lem:focal}, these are complementary spectral projections of
$A_\delta$. Hence
\[
 f(A_\delta)
 =
 f(\varsigma)P_++f(-\varsigma)P_-
 =
 a_fI+b_fA_\delta.
\]
If $\delta=0$, then $A_0^2=0$, and $f(A_0)=f(0)I+f'(0)A_0$. In both cases,
$W(f(A_\delta))=a_f+b_fK_\delta$.

\smallskip
\noindent\emph{Step 2: an affine test.} Fix $T$ with
$W(T)\subseteq K_\delta$, and put
\[
 N_f=\alpha T+\alpha' I,
 \qquad
 D_f=\gamma T+\gamma' I.
\]
Since $\sigma(T)\subseteq K_\delta$, the operator $D_f$ is invertible, and
$f(T)=N_fD_f^{-1}$. For a unit vector $x$, write $z=\inner{Tx}{x}$ and
$b=\norm{Tx}^2$. Then
\begin{align*}
 \norm{N_fx}^2
 &=
 \abs{\alpha}^2b+2\Real(\alpha\overline{\alpha'} z)+\abs{\alpha'}^2,\\
 \norm{D_fx}^2
 &=
 \abs{\gamma}^2b+2\Real(\gamma\overline{\gamma'} z)+\abs{\gamma'}^2,\\
 \inner{N_fx}{D_fx}
 &=
 \alpha\overline\gamma b+\alpha\overline{\gamma'} z
 +\alpha'\overline\gamma\overline z+\alpha'\overline{\gamma'}.
\end{align*}
Thus $\norm{N_fx}^2$, $\norm{D_fx}^2$ and the real and imaginary parts of
$\inner{N_fx}{D_fx}$ are real affine functions of $(\Real z,\Imag z,b)$.

By Proposition~\ref{prop:shell}, such a function is nonnegative on the
union of the shells as soon as it is nonnegative at the two endpoints
\[
 (z,\abs{z}^2),
 \qquad
 (z,\mathfrak b_\delta(z))
\]
of every fiber, that is, for the scalar $z$ and for $A_\delta$ at the
vector of Proposition~\ref{prop:shell}.

\smallskip
\noindent\emph{Step 3: the norm.} Let $h$ be the right-hand side
of~\eqref{eq:individual}. If $h=0$, then $f=0$ on $K_\delta$, and there is
nothing to prove. Otherwise, by Lemma~\ref{lem:pencil}{\rm(i)} and
congruence by $D_f$, one has $w_\rho(f(T))\leq h$ if and only if
\[
 \rho h^2\norm{D_fx}^2
 -2ph\Real\bigl(\overline\zeta\inner{N_fx}{D_fx}\bigr)
 -q\norm{N_fx}^2
 \geq0,
 \qquad
 \norm{x}=1,
 \quad
 \zeta\in\T.
\]
For the scalar $z$, this holds by Lemma~\ref{lem:pencil}{\rm(i)} applied
to the $1\times1$ operator $f(z)$, since $w_\rho(f(z))=\abs{f(z)}\leq h$.
For $A_\delta$, it holds because $w_\rho(f(A_\delta))\leq h$. By Step~2,
\[
 w_\rho(f(T))\leq h
 \qquad\text{whenever}\quad
 W(T)\subseteq K_\delta.
\]
Conversely, the scalars $z\in K_\delta$ and $A_\delta$ belong to this
class. Lemma~\ref{lem:closed-class} identifies the supremum over this
class with $\norm{f}_{\rho,\mathrm{bfd}}$.

\smallskip
\noindent\emph{Step 4: the numerical range.} Fix $\zeta\in\T$ and
$h\in\R$ such that the half-plane
\[
 \{u:\Real(\overline\zeta u)\leq h\}
\]
contains $f(K_\delta)$ and $W(f(A_\delta))$. By congruence with $D_f$, the
inequality $\Real(\overline\zeta f(T))\leq hI$ is equivalent to
\[
 h\norm{D_fx}^2-\Real\bigl(\overline\zeta\inner{N_fx}{D_fx}\bigr)\geq0,
 \qquad
 \norm{x}=1.
\]
This holds for the scalars and for $A_\delta$ by the choice of $h$, and
hence for all $T$ by Step~2. Intersecting these half-planes places
$W(f(T))$ in the right-hand side of~\eqref{eq:individual-range}.

Conversely, a point of that set is a finite convex combination
\[
 \sum_j\theta_jw_j,
 \qquad
 w_j=f(z_j)\ \ (z_j\in K_\delta)
 \quad\text{or}\quad
 w_j=\inner{f(A_\delta)y_j}{y_j}\ \ (\norm{y_j}=1).
\]
The direct sum of the corresponding scalars and copies of $A_\delta$ has
numerical range in $K_\delta$. The direct sum of the vectors
$\sqrt{\theta_j}$ and $\sqrt{\theta_j}\,y_j$ is a unit vector, at which the
numerical range of $f$ of this operator takes the value
$\sum_j\theta_jw_j$.
\end{proof}

\begin{remark}\label{rem:lf}
For each fixed $\eta$, Theorem~\ref{thm:lf} applies to $f_{g_\eta}$
in~\eqref{eq:mob-lift} and describes its range over all operators with
$W(T)\subseteq K_\delta$. In Theorem~\ref{thm:shift}, the optimal $\eta$
depends on the shift, and the upper bound of Section~\ref{sec:block} shows
that no other symbol in $\Sch$ does better.
\end{remark}

\subsection{Focal rigidity and equality}\label{ssec:rigidity}

We first record the equality case of the upper estimate in
Proposition~\ref{prop:shell}.

\begin{lemma}\label{lem:focal-shell-rigidity}
Let $W(T)\subseteq K_\delta$, and let $x$ be a unit vector with
\[
 \norm{Tx}^2=\mathfrak b_\delta\bigl(\inner{Tx}{x}\bigr).
\]
Then $\operatorname{span}\{x,Tx\}$ is a two-dimensional reducing
subspace for $T$, and the restriction of $T$ to this subspace is
unitarily equivalent to $A_\delta$.
\end{lemma}

\begin{proof}
By Theorem~\ref{thm:choi-li}, we may write
\[
 T=J^*BJ,
 \qquad
 B=A_\delta\otimes I_{\mathcal M},
\]
where $J$ is an isometry. Write $Jx=(\xi_1,\xi_2)$, and put
\[
 \vartheta=\norm{\xi_2}^2,
 \qquad
 \chi=2\inner{\xi_2}{\xi_1}.
\]
As in the proof of Proposition~\ref{prop:shell},
\[
 \norm{Tx}^2
 \leq
 \norm{BJx}^2
 =
 4\delta^2+4(1-\delta^2)\vartheta
 \leq
 \mathfrak b_\delta\bigl(\inner{Tx}{x}\bigr).
\]
By assumption, both inequalities are equalities.

Equality in the second inequality gives $\vartheta=\vartheta_+$, that is,
$\abs{\chi}^2=4\vartheta(1-\vartheta)$. Hence equality holds in the Cauchy--Schwarz
inequality
\[
 \abs{\inner{\xi_2}{\xi_1}}\leq\norm{\xi_2}\norm{\xi_1},
\]
and $\xi_1,\xi_2$ belong to a common one-dimensional subspace $\C e$ of
$\mathcal M$. Equality in the first inequality gives
$BJx\in J\mathcal H$, because $JJ^*$ is the orthogonal projection onto
$J\mathcal H$. Hence
\[
 BJx=JJ^*BJx=JTx.
\]

Put $z=\inner{Tx}{x}$. Since $\abs{z}\leq1+\delta$,
\[
 \mathfrak b_\delta(z)-\abs{z}^2
 \geq
 2(1+\delta^2)-(1+\delta)^2
 =
 (1-\delta)^2
 >0.
\]
Thus $Tx$ is not a multiple of $x$, and therefore $Jx$ and $BJx=JTx$ are
linearly independent. Both vectors belong to the two-dimensional space
$\C^2\otimes\C e$, so they span it.

This subspace reduces $B$, and it is contained in $J\mathcal H$.
Therefore it also reduces the compression of $B$ to $J\mathcal H$. Since
$J$ is a unitary map of $\mathcal H$ onto $J\mathcal H$ that intertwines
$T$ with this compression, the subspace
\[
 J^{-1}(\C^2\otimes\C e)=\operatorname{span}\{x,Tx\}
\]
reduces $T$, and $T$ restricted to it is unitarily equivalent to
$A_\delta$.
\end{proof}

The next result describes equality for a fixed linear-fractional
function. It applies in particular to $f=f_{\beta_s}$; see
Corollary~\ref{cor:canonical-equality}.

\begin{theorem}\label{thm:lf-equality}
Let $1\leq\rho\leq2$, and let
\[
 f(z)=\frac{\alpha z+\alpha'}{\gamma z+\gamma'}
\]
be nonconstant, with $\gamma z+\gamma'\neq0$ on $K_\delta$. Put
\[
 h_{\rm sc}=\max_{z\in K_\delta}\abs{f(z)},
 \qquad
 h_{\rm fo}=w_\rho(f(A_\delta)),
 \qquad
 h=\max\{h_{\rm sc},h_{\rm fo}\}.
\]
Let $T$ act on a finite-dimensional space and satisfy
$W(T)\subseteq K_\delta$.

\noindent{\rm (i)} If $h_{\rm sc}>h_{\rm fo}$, then $w_\rho(f(T))=h$ if
and only if $T$ has a reducing scalar summand $z_0I$ with
$z_0\in\partial K_\delta$ and $\abs{f(z_0)}=h$.

\noindent{\rm (ii)} If $h_{\rm fo}>h_{\rm sc}$, then $w_\rho(f(T))=h$ if
and only if $T$ has a reducing summand unitarily equivalent to
$A_\delta$.

\noindent{\rm (iii)} If $h_{\rm sc}=h_{\rm fo}$, then $w_\rho(f(T))=h$ if
and only if $T$ has one of the two summands in~{\rm (i)} and~{\rm (ii)},
or there are $z_0\in W(T)$ and $\zeta\in\T$ such that
\[
 \abs{f(z_0)}=h,
 \qquad
 p\bigl(h-\Real(\overline\zeta f(z_0))\bigr)=0,
\]
and
\[
 \rho h^2\abs{\gamma}^2
 -2ph\Real(\overline\zeta\alpha\overline\gamma)
 -q\abs{\alpha}^2
 =0.
\]
\end{theorem}

\begin{proof}
Put $X=f(T)$, $N=\alpha T+\alpha'I$ and $D=\gamma T+\gamma'I$.
The denominator $D$ is invertible and $XD=N$. By~\eqref{eq:individual},
$w_\rho(X)\leq h$, and $h>0$ since $f$ is nonconstant. Write
\[
 Q_Y(\zeta,a)=\rho a^2I-2pa\Real(\overline\zeta Y)-qY^*Y .
\]
In finite dimension, equality $w_\rho(X)=h$ holds precisely when
$Q_X(\zeta,h)$ has a kernel for some $\zeta\in\T$.
Indeed, if all these pencils are positive definite, compactness gives a
uniform positive lower bound, so $h$ can be decreased in~\eqref{eq:pencil}.
Conversely, whenever $w_\rho(Y)\leq a<h$ with $a>0$, we have
$\Real(\overline\zeta Y)\leq aI$ and
\[
 Q_Y(\zeta,h)-Q_Y(\zeta,a)
 \geq (h-a)(\rho h+qa)I.
\]
Thus every pencil at $h$ is strictly positive in this case.

For a unit vector $x$, put $z=\inner{Tx}{x}$ and $b=\norm{Tx}^2$.
Congruence by $D$ gives
\begin{align*}
 P_\zeta(z,b)
 &:=
 \inner{Q_X(\zeta,h)Dx}{Dx}\\
 &=
 \rho h^2\norm{Dx}^2-2ph\Real(\overline\zeta\inner{Nx}{Dx})
 -q\norm{Nx}^2.
\end{align*}
The squared norms and the real and imaginary parts of $\inner{Nx}{Dx}$
are real affine functions of $(\Real z,\Imag z,b)$. For instance,
\[
 \inner{Nx}{Dx}
 =
 \alpha\overline\gamma b+\alpha\overline{\gamma'}z
 +\alpha'\overline\gamma\,\overline z+\alpha'\overline{\gamma'}.
\]
It follows that the coefficient of $b$ in $P_\zeta$ is
\[
 L_\zeta=\rho h^2\abs{\gamma}^2
 -2ph\Real(\overline\zeta\alpha\overline\gamma)-q\abs{\alpha}^2.
\]
Since $D$ is invertible, the kernel criterion is equivalent to
$P_\zeta(z,b)=0$ at the shell pair of some unit vector $x$.

Both endpoint values of $P_\zeta$ on the fiber over $z$ are
nonnegative by Proposition~\ref{prop:shell} and~\eqref{eq:pencil}.
At the lower endpoint, division by $\abs{\gamma z+\gamma'}^2$ gives
\[
 (h-\abs{f(z)})(\rho h+q\abs{f(z)})
 +2ph\bigl(\abs{f(z)}-\Real(\overline\zeta f(z))\bigr).
\]
Both terms are nonnegative, so this endpoint is a zero precisely when
\begin{equation}\label{eq:scalar-zero}
 \abs{f(z)}=h,\qquad p\bigl(h-\Real(\overline\zeta f(z))\bigr)=0.
\end{equation}
At the upper endpoint the pencil is evaluated on $A_\delta$. Its value
is strictly positive if $h>h_{\rm fo}$, by the strict inequality for
$Q_Y$ above applied to $Y=f(A_\delta)$.

Suppose first that $h_{\rm sc}>h_{\rm fo}$. Affinity in $b$ forces
every zero to lie at the lower endpoint. Hence $Tx=zx$ by equality in
Cauchy--Schwarz, and $\abs{f(z)}=h$. The maximum modulus principle
gives $z\in\partial K_\delta$. For an outward unit normal $\omega$ at
$z$, the operator $\Real(\overline\omega(zI-T))$ is positive and has
zero quadratic form at $x$. It therefore annihilates $x$. Using
$Tx=zx$ gives $T^*x=\overline z x$, as required in~{\rm(i)}.
If $h_{\rm fo}>h_{\rm sc}$, the lower endpoint is strictly positive
on every fiber. A zero must then lie at the upper endpoint, and the
rigidity lemma~\ref{lem:focal-shell-rigidity} gives~{\rm(ii)}.
Conversely, either extremizing reducing summand gives equality because
the functional calculus respects direct sums and $w_\rho$ takes their
maximum.

Finally, suppose $h_{\rm sc}=h_{\rm fo}$. A zero at an endpoint gives
the corresponding summand by the same arguments. A zero strictly
between the endpoints forces the entire affine function on that fiber
to vanish. This is equivalent to~\eqref{eq:scalar-zero} and
$L_\zeta=0$. Conversely, if these conditions hold at $z_0\in W(T)$,
choose a unit vector $x$ with $\inner{Tx}{x}=z_0$. Then
$P_\zeta(z_0,\norm{Tx}^2)=0$, which proves~{\rm(iii)}.
\end{proof}

\begin{example}
The third alternative at a tie can occur without either reducing
summand. Take $f(z)=z$, $\rho=2$ and
\[
 T=\theta(1+\delta)I+(1-\theta)A_\delta,\qquad 0<\theta<1.
\]
By convexity, $W(T)=\theta(1+\delta)+(1-\theta)K_\delta\subseteq K_\delta$,
and $1+\delta\in W(T)$. Thus $w(T)=h_{\rm sc}=h_{\rm fo}=1+\delta$.
Here $p=1$, $q=0$ and $\gamma=0$, so $L_1=0$ and
$z_0=1+\delta$ satisfies~\eqref{eq:scalar-zero}. The matrix $T$ is
nonnormal and hence has no reducing line. Its nonzero trace also
excludes unitary equivalence to $A_\delta$.
\end{example}

Finite dimension is essential for the reducing-summand conclusions.
For $\delta=0$, $f(z)=z$ and $\rho=1$, the operator
\[
 T=\bigoplus_{n\geq2}(1-1/n)A_0
\]
has $w(T)\leq1$ and $\norm{T}=2>h_{\rm sc}=1$. Every nonzero vector has
a nonzero component in a summand of norm strictly less than $2$, so no
unit vector attains $\norm{T}$. In particular $T$ has no reducing copy of
$A_0$. This is why shell contact in Theorem~\ref{thm:equality-contact}
is formulated using the closure.

\begin{corollary}\label{cor:canonical-equality}
Let $0\leq s<1$. For the fixed transform $f=f_{\beta_s}$, the two
values in Theorem~\ref{thm:lf-equality} are
\[
 h_{\rm sc}=R-C,
 \qquad
 h_{\rm fo}=\frac{r+p\abs{c}+\sqrt{(r+p\abs{c})^2-\rho q(r^2-d^2-c^2)}}{\rho}.
\]
In finite dimension, their strict ordering determines the reducing
summand required for equality.
\end{corollary}

\begin{proof}
Let $\lambda=u+\delta\overline u$ with $u\in\T$. A direct computation
shows that the disk automorphism $\beta_s$ maps the disk
$\{\abs{w}\leq\delta\}$ onto the closed disk with center
$s(1-\delta^2)/(1-k\delta^2)$ and radius $R-1$. Since
$s(1-\delta^2)/(1-k\delta^2)-s=C$ by~\eqref{eq:parameters}, it follows that
\[
 \beta_s(\delta\overline u)-s\in\overline D(C,R-1).
\]
As $\abs{\beta_s(u)}=1$, the boundary formula~\eqref{eq:faber-intro}
gives
\[
 \abs{f_{g_1}(\lambda)-C}
 =
 \abs{\beta_s(u)+\beta_s(\delta\overline u)-s-C}
 \leq
 1+(R-1)
 =
 R.
\]
Thus $h_{\rm sc}\leq R+\abs{C}=R-C$, with equality at
$\lambda=-(1+\delta)$ by Lemma~\ref{lem:scalar-model}. The focal value
$h_{\rm fo}=w_\rho(f_{g_1}(A_\delta))$ is the second value
in~\eqref{eq:radius}, by the proof of Corollary~\ref{cor:radius}.
\end{proof}

\begin{remark}\label{rem:not-lf}
For $\delta=0$ and $g(u)=u^2$, consider the operator of
Crabb~\cite{Crabb} that attains $\norm{T^2}=2$ under $w(T)\leq1$; see
also Gau and Wu~\cite[Corollary~4]{GauWu}:
\[
 T=\begin{pmatrix}0&\sqrt2&0\\0&0&\sqrt2\\0&0&0\end{pmatrix}.
\]
This weighted shift has $w(T)=\sqrt2\cos(\pi/4)=1$
(see~\cite{HaagerupDeLaHarpe}). Since $T^2e_3=2e_1$ and $T^2$ vanishes
on $\operatorname{span}\{e_1,e_2\}$, we have
\[
 \norm{T^2}=2=\mathfrak M_{0,0}(0),\qquad
 (\inner{T^2e_3}{e_3},\norm{T^2e_3}^2)=(0,4).
\]
Thus $g$ attains the strict focal value for the whole symbol class.
A reducing copy of $A_0$ would leave a reducing line in its orthogonal
complement. Nilpotence would place that line in
$\ker T\cap\ker T^*=\{0\}$.

For the fixed function $f(z)=z^2$, however, the maximum over scalars
is $1$ and $f(A_0)=0$. The same example therefore shows that the
fixed-function formula~\eqref{eq:individual} itself can fail for a
quadratic. Its norm involves the fourth-order moment $\norm{T^2x}^2$,
whereas congruence for a linear-fractional function uses only
$\inner{Tx}{x}$ and $\norm{Tx}^2$.
\end{remark}

\subsection*{Declarations}
\begin{itemize}
\item \textbf{Funding.} The author received no funding for this work.
\item \textbf{Competing interests.} The author declares that there are no
competing interests.
\item \textbf{Data availability.} Data sharing is not applicable to this
article, as no datasets were generated or analysed.
\item \textbf{Use of AI tools.} Claude (Anthropic) was used throughout
this work, under the author's direction, to develop and test arguments,
carry out calculations, and draft and revise the manuscript. The author
checked every proof and calculation line by line and takes full
responsibility for the content of this article.
\end{itemize}

\end{document}